\documentclass[fleqn,11pt]{article}
\usepackage{graphicx}
\usepackage{amscd,amssymb}
\usepackage{amsmath}
\title{Characters for Complex Bundles and their Connections}
\author{James Simons\thanks{Stony Brook University, Simons Foundation} and Dennis Sullivan\thanks{CUNY Graduate Center, Stony Brook University}}
\date{}

\begin{document}
\setcounter{footnote}{1}

\maketitle

\begin{abstract}
Given a complex vector bundle $E$ over a base manifold $X$ with connection $\nabla$ we construct invariants called characters in terms of integrals over manifolds with boundary:

\vspace{.5cm}

\begin{center}
\begin{tabular}{|c|c|}
\hline
$\nabla$ & character values \\
\hline
$\mathbb{C}$-linear connection & $\mathbb{C}/\mathbb{Z}$ \\
\hline
unitary connection & $\mathbb{R}/\mathbb{Z}$ \\
\hline
independent of connection & $\mathbb{Q}/\mathbb{Z}$ \\
\hline
\end{tabular}
\end{center}

\vspace{.5cm}

In logical order: the first result, Theorem AT (Algebraic Topology), shows the $\mathbb{Q}/\mathbb{Z}$-characters derived from $\mathbb{C}/\mathbb{Z}$-characters are in a bijective correspondence with complex $K$-theory.  The second result, Theorem DG (Differential Geometry), shows the $\mathbb{C}/\mathbb{Z}$ (or $\mathbb{R}/\mathbb{Z}$) characters are in a bijective correspondence with differential $K$-theory defined using complex (or real) valued differential forms.  Differential $K$-theory may be defined as the Grothendieck group of Chern-Simons equivalence classes of complex bundles with connection (respectively $\mathbb{C}$-linear or unitary).  The third result, Theorem AN (Analysis): \textit{i}) Expresses the unitary bijection in terms of the eta invariants mod one of the $\textrm{spin}^{c}$ Dirac operators with coefficients in $(E,\nabla)$ restricted to enriched closed odd-dimensional stably almost complex manifolds (SACs) in $X$, and \textit{ii})a) There is an easy and natural push forward in differential $K$-theory ($\mathbb{C}$-linear or unitary) for families of even-dimensional SACs defined by pulling back the odd-dimensional SAC cycles in $X$ which are enriched using the direct sum connection, \textit{ii})b) in the unitary case the direct sum connection is $CS$ equivalent to a rescaling (adiabatic) limit of the Levi-Civita connection, and one thereby computes the push forward (unitary case) as a limit of eta invariants mod one as the rescaling tends to infinity.

This adds to the discussion in the literature of the question (communicated by Iz Singer) to have an index theorem in differential $K$-theory for families.

\vspace{.5cm}
History and Acknowledgement:  The new aspect of Theorem AT is the construction and proof of the bijection using integrals of characteristic forms, whereas the existence of such a bijection is a variant of much earlier understanding [1].  The form of the $\hat{K}$-character definition is a modified form of ``Differential Characters'' [2] motivated by ``The Characteristic Variety Theorem'' [3].  We are indebted to [4] for completing [5] for the non-unitary, $\mathbb{C}$-linear connections describing differential $K$-theory with $\mathbb{C}$-valued differential forms.
\end{abstract}

\section*{Introduction}

%\renewcommand{\thefootnote}{\fnsymbol{footnote}}
%\footnotetext[1]{Jim's address}
%\footnotetext[2]{Dennis' address}
\renewcommand{\thefootnote}{\arabic{footnote}}

Let $V$ be an odd-dimensional Stably Almost Complex (SAC) closed
manifold mapping smoothly $F : V \rightarrow X$ to a target manifold
$X$.  Suppose $X$ is the base of a complex vector bundle $E$ together with a $\mathbb{C}$-linear connection $\nabla$.  Suppose $V$ is ``enriched'' by an independent $\mathbb{C}$-linear connection $\nabla/V$ on its stable tangent bundle. We can form the Chern character, a total even form $ch(E)$ on $X$.  We can also form on $V$ the total even Todd form of $V$ (see background in Section 0).

One can show [Section 1] that $V$ bounds in the SAC sense a $W$ in
such a way that the bundle $F^{\ast}(E)$ on $V$ extends to a complex
vector bundle  $E/W$ on $W$.  One can impose a $\mathbb{C}$-linear
connection on the stabilized tangent bundle of $W$ compatible with
that on $V$ in a collar neighborhood of $V$, and also a connection on
$E/W$ similarly compatible with that on $F^{\ast}(V)$.  Using these
connections we can extend $\textrm{Todd}(V)$ and $F^{\ast}(ch(E))$ to total
forms $\textrm{Todd}(W)$ and $ch(E/W)$ over all of $W$.
Now form the integral over $W$ of the wedge product $\textrm{Todd}(W)
ch(E/W)$ and reduce the value modulo one.  We refer to this as an
``angle'' in the complexified circle $\mathbb{C}/\mathbb{Z}$.  This ``angle'' in $\mathbb{C}/\mathbb{Z}$ may be seen to be independent of the choice
of $W$ and the extended connections above.   It depends only on $F : V
\rightarrow X$ and the enrichment of $V$.  As the triple $({V,F,\nabla/V})$ varies one obtains many ``angle'' invariants of $(E,\nabla)$, the complex bundle with $\mathbb{C}$-linear connection over $X$. The ``angle'' invariants are thought of as a ``character'' function with values in the complexified circle $\mathbb{C}/\mathbb{Z}$ defined on the collection of odd-dimensional enriched SAC closed manifolds or cycles mapping into $X$.

One sees directly the ``angle'' invariants of $(E,\nabla)$ over $X$ are equal as ``characters'' to the ``angle'' invariants of $(E,\nabla^{\prime})$ if $\nabla$ and $\nabla^{\prime}$ are $\mathit{CS}$ equivalent.  This means the odd-dimensional Chern Simons forms whose exterior $d$ is the difference between the corresponding Chern character forms are not only closed so that the Chern character forms of $(E,\nabla)$ and $(E,\nabla^{\prime})$ are exactly the same, but these $CS$ forms are also exact. This statement of equality of ``angle invariants'' for $CS$ equivalent connections follows from Stokes' theorem.

One may also change the pair $(E,\nabla)$ by direct summing with the trivial connection on the trivial bundle or by changing by a strict isomorphism of bundles with connection, without changing the values of the character function associated with $(E,\nabla)$.

The goal of this paper is twofold.  Firstly, to prove the converse of these equalities of ``angle'' invariants.  Namely, the ``angle'' invariant characters actually determine $(E,\nabla)$ up to the equivalences just mentioned.  Secondly, to describe the properties of the character function on enriched SAC cycles in $X$ that are necessary and sufficient for such a function to arise as the ``angle'' invariants of a complex bundle $E$ with $\mathbb{C}$-linear connection over $X$. These properties are described / contained in the differential geometry statement \textbf{Theorem DG} in Section 3.

First goal: the form of the equivalences was exploited recently [5] to give a geometric description, in the case of the generalized cohomology theory complex $K$-theory, of the extension introduced in [6] of a contravariant functor called ``generalized differential cohomology''.  This is a derived fibre product of the generalized cohomology functor $h^{\ast}$ with de Rham differential forms either real or complex valued which are labeled by elements in $h^{\ast}(pt)$. The fibre product is over the generalized de Rham isomorphism canonically identifying $h^{\ast}(\; ,\mathbb{R} \; \textrm{or} \; \mathbb{C})$ with the de Rham cohomology of forms (over $\mathbb{R}$ or $\mathbb{C}$) with coefficients $h^{\ast}(pt)$.  The geometric description of [5] says  differential complex $K$-theory defined by forms with values in $\mathbb{R}$ or $\mathbb{C}$, denoted $\hat{K}$, is represented by the Grothendieck group of stable isomorphism classes of complex bundles with $\mathbb{C}$-linear connections in the case of $\mathbb{C}$-valued forms and with unitary connections in the case of $\mathbb{R}$-valued forms, up to $\mathit{CS}$ equivalence (structured bundles) under direct sum. The proof of the odd form lemma in [5] depended on the existence of inverses up to $CS$ equivalance for the bundles with connection.  These inverses were provided for unitary connections in [5].  For $\mathbb{C}$-linear connections, the inverses were provided in [4].  We are indebted to Leon Takhtajan for explaining how [4] completes the $\mathbb{C}$-linear results of [5].  In particular, Corollary 2 of [4] plus Proposition 2.5 of [5] implies Corollary 3 of [4] which produces the required inverses.  There is yet another proof for GL[N,C] of the required inverse property of bundles with connection in [7]. 

A helpful organizing tool (Section 3) in the proof of the differential geometry statement Theorem DG is the ``character'' diagram of groups placing $\hat{K}$ in the hexagons between the interlocking sine wave Bockstein long exact sequence and the de Rham long exact sequence.  Verifying the analogous diagram for $\hat{K}$ characters is the rest of the work in the proof of Theorem DG beyond Theorem AT mentioned above and described more fully below.

Second goal: the pattern of the necessary and sufficient conditions on
``angle'' invariants to come from a bundle with $\mathbb{C}$-linear
connection modifies the notion of ``differential character''
introduced in [2]. These were functions $C$ with values in
$\mathbb{R}/\mathbb{Z}$ on smooth cycles in $X$ in the sense of
algebraic topology. These functions were additive under union of
cycles. They were not homology invariant but their values in
$\mathbb{R}/\mathbb{Z}$ changed under a homological deformation of a
cycle by integrating a closed real differential form with integral
periods on $X$ over the homology.  Namely if $z$ and $z^{\prime}$ are
cycles in $X$ then $C(z) - C(z^{\prime})$ in $\mathbb{R}/\mathbb{Z}$
is given by the integral of $w(C)$ over any homology in $X$ between
$z$ and $z^{\prime}$ reduced [mod one]. In particular a small smooth
deformation produces a small change in the value, thus the name
differential character.  This definition was motivated by the attempt
to define objects in the base related to the Chern-Simons forms in the bundle.

The modification of the notion of differential character to the main
notion of differential $\hat{K}$ character employed here to describe
the ``angle'' invariants goes as follows: replace algebraic topology
cycle in $X$ as in [2] by smooth enriched odd-dimensional SAC cycle in
$X$, and replace the variation $k$-form $w(C)$ of the differential
character by a total even complex valued form $W(C)$ representing the
Chern character of some complex vector bundle over $X$.  Replace the
variation formula for differential characters by $C(V) -C(V^{\prime})$
equals the integral mod $\mathbb{Z}$ of  $W(C) \wedge \textrm{Todd} \;
M$ over $M$, where $M$ is an enriched SAC mapped into $X$ whose boundary is the formal difference of enriched cycles $V - V^{\prime}$ in $X$.  A new ``product relation'' beyond [2] appears because of the multiplicative properties of the Todd form.  This modification follows the pattern of the Characteristic Variety Theorem [3].

Theorem DG gives a second geometric interpretation added to that of [5] for differential $K$-theory of $X$ in terms of differential $\hat{K}$ characters defined on enriched SAC cycles in $X$.  In the real case using unitary connections, thanks to the APS theorem [8], this bijection is elegantly described by forming the eta invariants mod one of the $\textrm{spin}^{c}$ Dirac operator of the odd-dimensional SAC cycle with coefficients in the bundle $E$ restricted to $V$ to build the $\hat{K}$ character of the bundle with unitary connection.

\pagebreak

The proof scheme for Theorem DG uses the machinery of algebraic topology to arrive first at a corresponding result, Theorem AT (Algebraic Topology), giving a complete theory of rational ``angles'' associated to just the bundle $E$ independent of the choice of connection.  The rational ``angles'' of order $k$ in $\mathbb{C}/\mathbb{Z}$ are associated with pairs $(V,W)$ where $V$ is an odd SAC cycle in $X$ with a given way, $W$, to bound $k$ copies of $V$ in $X$.  These rational angles are derived in Section 2, in terms of the general ``angles'' associated with enriched SAC cycles in $X$ probing bundles with connections.  As the enrichment connections vary continuously, the rational ``angles'' cannot vary continuously without staying constant and thus become topological invariants of the bundle. They are complete invariants and their precise necessary and sufficient conditions are specified using Rational $\textrm{Hom}(h( \; ,\mathbb{Q}/\mathbb{Z}), \mathbb{Q}/\mathbb{Z})$.  Here $h$ is the remarkable homology theory discovered by Conner and Floyd [9] using SAC cycles and homologies (i.e. SAC bordism classes) in $X$ taken also modulo the purely algebraic ``product relation'' mentioned above. The fact that this algebraic quotient does not destroy the exactness property of a homology theory is the crucial point of the argumentation here.  This is the Algebraic Topology theorem, Theorem AT, described in Section 2 and used in Section 3.  

One interesting corollary of the Topology discussion is 

\underline{\textbf{Corollary}}:  A cohomology class $c$ in $H^{even}(X, \mathbb{Q})$ is the Chern character of a complex bundle over $X$ if and only for every closed even-dimensional SAC mapping to $X$, $V \xrightarrow{f} X$, $\int f^{\ast} \; c \; \textrm{Todd} \; V$ is an integer.  A similar statement holds for the transgressed $ch$ in $U$, odd-dimensional closed SACs in $X$, elements in $H^{odd}(X, \mathbb{Q})$ and maps $X \rightarrow U$.

The first application of Theorem DG is a very easy explicit construction of a wrong way map in differential $K_{\mathbb{C}}$-theory for a fibration with $\mathbb{C}$-linear connection over a base $X$ with fibres closed even-dimensional SAC manifolds enriched by $\mathbb{C}$-linear connections.  If $T$ denotes the total space there is a map of enriched SAC cycles from $X$ to those in $T$ by taking the pullback SAC cycle.  There is a nuance here, but then applying $\hat{K}$ characters  reverses the direction to give a wrong way map from the differential $K$-theory of $T$ to the differential $K$-theory of $X$ (all over $\mathbb{C}$).  See [14].  The nuance here is the enrichment of the pullback cycle.  Since the fibration is enriched with a connection on the vertical subbundle, one may use the direct sum connection (see Section 4) to enrich the pullback cycles. Then a differential $\hat{K}$ character on $T$ will induce a differential $\hat{K}$ character on the base using the multiplicative nature of Todd forms and Chern character forms to define the required variation form.  (See Section 4.)  

There is also the statement:

\underline{\textbf{Corollary 1 of Proof of Theorem DG}}: \linebreak kernel$(\hat{K}X \xrightarrow{ch} \wedge^{even}_{integrality})$ is isomorphic to $\textrm{Hom}(\bar{\Omega}_{odd}^{\mathbb{C}},\mathbb{C}/\mathbb{Z})$, a complex torus of dimension the sum of the odd Betti numbers of $X$.

The final discussion brings in analysis and the eta invariants of the Atiyah-Patodi-Singer [8] theory relating Topology, Geometry, and Analysis of a SAC manifold with boundary, but now restricted to using unitary connections on the bundles and Levi-Civita connections on manifolds.

\pagebreak

As mentioned above, the ``angle invariant'' in $\mathbb{R}/\mathbb{Z}$ for a unitary connection constructed using the filling $W$ is just the spectral invariant on $V$ defined using the $\textrm{spin}^{c}$ Dirac operator reduced modulo one. Here it is important though that the connection on $V$ be the Levi-Civita connection so that the asymptotic heat kernel analysis of the APS Theorem is valid.

Fortunately but not obviously, there is an extension of this
calculation relating real ``angles'' to eta invariants mod one giving
an analytic calculation of the wrong way map in differential
$K$-theory in the case of unitary connections.  There is a serious
stumbling block though.  For the APS theory one again needs the
asymptotic analysis based on using the Levi-Civita connection on the
total space $T$.  For the wrong way map one needs the multiplicative
property of the  characteristic forms associated to the direct sum
connection on $T$ in order to define the variation form in the
definition of $\hat{K}$-character.  Even though the metric on $T$ can
be taken to be the direct sum metric the Levi-Civita connection is not
the direct sum connection.  See the Appendix to Section 5 for the
detailed discussion of the interesting difference.  However in the
Appendix one sees that by scaling in such a way that the base becomes
infinitely large compared to the fibre the limit of the Levi-Civita
connections upstairs exists (the adiabatic limit)\footnote{We learned
  of the adiabatic limit connection from the work of Cheeger [12] and
  Freed [17].} and is fortunately $CS$ equivalent\footnote{The
  equality of the Chern character forms was known before
  [communication of Dan Freed] but $CS$ equivalence is a stronger condition.} to the direct sum connection. This
equivalence enables the analytic calculation of the wrong way map in differential $K$-theory (unitary case) as a limit of eta invariants of the rescaled metrics. This is the analytic theorem, Theorem AN, in Section 5.

Theorem AN is our response to Iz Singer's question (on a flight with the authors seven years ago) about having an analytic version of the index theorem in differential $K$-theory for families. There have been other responses which use infinite-dimensional analysis [10] and [references in [11]].

\section*{\S0. Background about connections, Chern-Weil characteristic forms, and Chern-Simons forms}

Let $W \xrightarrow{\Pi} X$ be a real $n$-dim vector bundle over a
smooth manifold.  Let $\nabla$ be a connection, and $R \in
\wedge^{2}(X, \mathrm{End}(W))$ the curvature tensor.

A real valued polynomial on the Lie algebra of $GL(n,R)$ is called
\textbf{invariant} if it is fixed under the adjoint action.  If
$P_{l}(B) = t_{r}(B^{l})$, it is well known that the ring of invariant
polynomials is generated by $P_{1}, \cdots, P_{n}$.  For a degree $l$
invariant polynomial $P$, the Chern-Weil homomorphism yields
\begin{eqnarray*}
P(\overbrace{R \wedge R \wedge \cdots \wedge R}^{l}) \in
\wedge^{2l}(X), \qquad \textrm{a closed form.}
\end{eqnarray*}

This map is a ring homomorphism, and the associated cohomology class
of an element in its image is independent of the choice of
connection.  This is made apparent by 1.1) and 1.2) below.

If $\nabla^{0}$ and $\nabla^{1}$ are two connections on $W$ with
curvature forms $R^{0}$ and $R^{1}$, and $\nabla^{t}$ is a smooth curve
of connections joining $\nabla^{0}$ to $\nabla^{1}$, and $R^{t}$ its
curvature, set
\begin{eqnarray*}
B^{t} = \frac{d}{dt}(\nabla^{t}) \in \wedge^{1}(X,\mathrm{End}(W)).
\end{eqnarray*}
For $P$ invariant of deg $l$,
\[
1.1) \qquad P(R^{1} \wedge \cdots \wedge R^{1}) - P(R^{0} \wedge \cdots \wedge R^{0}) = d(\mathit{TP}(\nabla^{0},\nabla^{1})) \qquad \mathrm{where}
\]
\[
1.2) \qquad \mathit{TP}(\nabla^{0},\nabla^{1}) = l \int_{0}^{1} P(B^{t} \wedge
\overbrace{R^{t} \wedge \cdots \wedge R^{t}}^{l-1}) \; \bmod \wedge^{2l-1}_{\textrm{exact}}. \hspace{.3cm} \textrm{(Chern-Simons forms)}
\]
It may be shown that $\mathit{TP}$ is independent of the curve joining
$\nabla^{0}$ to $\nabla^{1}$, and thus $\mathit{TP}(\nabla^{0},\nabla^{1})$ is
well defined.  We also have
\[
1.3) \qquad \mathit{PQ}(R \wedge \cdots \wedge R) = P(R \wedge \cdots \wedge R)
\wedge Q(R \wedge \cdots \wedge R)
\]
\[
1.4) \qquad T(\mathit{PQ})(\nabla^{0},\nabla^{1}) = \mathit{TP}(\nabla^{0},\nabla^{1})
\wedge Q(R^{0} \wedge \cdots \wedge R^{0}) + \mathit{TQ}(\nabla^{0},\nabla^{1})
\wedge P(R^{1} \wedge \cdots \wedge R^{1}).
\]
The first, because the Chern-Weil map is a ring homomorphism, and the
second by calculation (recall $\mathit{TP}$ is defined $\bmod$ exact).

\underline{\textbf{Definition}}:  $\nabla^0$ and $\nabla^1$ are called $\mathbf{CS}$
\textbf{equivalent} if $\mathit{TP}(\nabla^0, \nabla^1)$ is exact for all
invariant $P$.  This is easily shown to be an equivalence relation.

Since $\{P_{l}\}$ generates the ring of invariant polynomials, 1.3)
and 1.4) show

\underline{\textbf{Proposition 1.5}}:  $\nabla^0$ and $\nabla^1$ are
equivalent if and only if $\mathit{TP_{l}}(\nabla^0, \nabla^1)$ is exact for
all $l \leq n$.

\underline{\textbf{Remark}}:  The discussion for $\mathbb{C}$-linear connections and complex-valued Chern-Weil and Chern-Simons is similar.

\section*{\S1. Construction I for Stably Almost Complex Manifolds (SACs)}

\underline{\textbf{Proposition 1}}:  Any complex vector bundle $E$ over $\Sigma$, a closed odd-dimensional SAC, can be filled in.  Namely, there is an even-dimensional SAC, $W$ with $\partial W = \Sigma$, and a complex bundle $E_{W}$ over $W$ extending $E$.

\underline{\textbf{Proof}}:  The SAC bordism of a point, $\Omega^{\mathbb{C}}_{\ast}(pt)$ is torsion free and is concentrated in even degrees [a celebrated result of the 60's]. So is the homology of $B_{U_{n}}$ the classifying space of isomorphism classes of complex bundles of rank $n$.  The SAC bordism of $B_{U_{n}}$, $\Omega^{\mathbb{C}}_{\ast}(BU_{n})$, is the limit of the Atiyah-Hirzebruch spectral sequence which begins with $H_{\ast}(BU_{n},\Omega^{\mathbb{C}}_{\ast}(pt))$ and which for any homology theory collapses when tensored with $Q$.  Thus it already collapses in this case and $\Omega^{\mathbb{C}}_{\ast}BU_{n}$ is torsion free and concentrated in even degrees.  This proves Proposition 1.

Construction I below defines a pairing $<\nabla_{\Sigma}; \Sigma \xrightarrow{f} X | E \rightarrow X; \nabla_{E}>$ in $\mathbb{C}/\mathbb{Z}$ where $\Sigma \xrightarrow{f} X$ is a closed odd-dimensional SAC in $X$, $E \rightarrow X$ is a complex vector bundle over $X$, $\nabla_{\Sigma}$ is a $\mathbb{C}$-linear connection on the stable tangent bundle of $\Sigma$, and $\nabla_{E}$ is a $\mathbb{C}$-linear connection on $E$ over $X$.

\underline{\textbf{Construction I}}:  By Proposition 1 we can fill in $\Sigma$ by $W$ and extend $E_{\Sigma} = E / \Sigma$ to $E_{W}$ over $W$.  Similarly we can extend $\nabla_{\Sigma}$ to $\nabla_{W}$ on the stable tangent bundle of $W$ and $\nabla_{E} / \Sigma$ to the extended bundle $E_{W}$ over ${W}$. We suppose there is a product neighborhood near $\partial W$ where the extended connections are product-like.  On $W$ there are two even-dimensional differential forms $ch(E_{W},\nabla)$, the Chern character form defined by $(E_{W},\nabla)$, and Todd $W$, the characteristic form associated to the universal total Todd class which is constructed from $\nabla_{W}$. (To be precise, we use the Todd form, as defined in the Remark below, associated with the inverse of the stable tangent bundle.)

Define the pairing 
\[
<\nabla_{\Sigma}; \Sigma \xrightarrow{f} X | E \rightarrow X; \nabla_{E}> \, \mathrm{in} \, \mathbb{C}/\mathbb{Z}
\]
by $\int_{W} ch(E_{W}) \cdot \mathrm{Todd} W$ reduced modulo 1.

\underline{\textbf{Proposition 2}}:  The value of the integral of $W$ mod one only depends on the SAC cycle in $X$, $\Sigma \xrightarrow{f} X$ ``enriched'' by the connection $\nabla_{\Sigma}$ and the complex vector bundle $E \rightarrow X$ ``enriched'' by its connection $\nabla_{E}$.

\underline{\textbf{Proof}}:  If we had chosen a different filling $\bar{W}, E_{\bar{W}}$ so that $\partial \bar{W} = \Sigma$ and $E_{\bar{W}} / \Sigma = E_{\Sigma}$ and different ``enrichment'' $\nabla_{\bar{W}}$ and $\nabla E_{\bar{W}}$ extending the enrichments on $\Sigma$ we can form the union of these two choices along $\Sigma$, namely $W \cup_{\Sigma} \bar{W}$ and $E_{W} \cup_{E_{\Sigma}} E_{\bar{W}}$ and also glue the enrichments.  The difference of the two integrals which are the two definitions of the pairing is the entire integral over the closed manifold $W \cup_{\Sigma} \bar{W}$.  This integral is well known to be an integer.  (See the next remark for some history.)  This proves Proposition 2.

\underline{\textbf{Remark 1}}:  Here we fix the definition of the Todd
class and recall how the integrality of $\int_{V} ch(E) \mathrm{Todd}
V$ was understood for closed even-dimensional SACs $V$.  

Imagine $V$ embedded in a large sphere $S$ with a complex structure on
its normal bundle $\mathcal{V}$ provided with a unitary structure.
Pull back $\mathcal{V}$ to the normal disk bundle $N$ to obtain a
complex bundle $E$ on $N$.  For each point $v$ of $N$ not in the zero
section of $\nu$ there is an isomorphism between the two halves of the
exterior algebra bundle associated to $E$, $\wedge^{even} E
\leftrightarrow \wedge^{odd} E$ defined by (wedging with $v$)
plus (contracting with $v$).

This defines an element in the complex $K$-theory of the pair
$K^{even}_{\mathbb{C}}$(disk bundle, sphere bundle) which we
pull back to the big sphere $S$ by the collapsing map $S
\leftrightarrow $ disk bundle / sphere bundle.  

The Chern character of this pullback element in the sphere is an
integer [Bott, Milnor, Adams].  One calculates in the universal
example over $B_{U}$ that $ch(\wedge^{even} - \wedge^{odd})$ in the
cohomology of the universal Thom space $MU$ satisfies
\[
ch(\wedge^{even} - \wedge^{odd}) = U \cdot \mathrm{Todd},
\]
where $U$ is the Thom class, for some universal class Todd (defined by this equation).  This
reveals the integrality above.

Previously, in his celebrated treatise, Hirzebruch calculated Todd in
terms of the multiplicative series $x/e^{x}-1$, motivated by Todd's
work on the arithmetic genus (or holomorphic Euler characteristics of
algebraic varieties).  This series shows the multiplicative property of the Todd formula.

The integrality made precise by Adams, Bott, and Milnor circa 1960
inspired Atiyah and Singer to build Dirac operators and to develop the
index theorem. (Recounted to one of the authors by Iz Singer (late 60's))  We revisit this later.

\underline{\textbf{Remark 2}}:  There are several properties of the
pairing $<\nabla_{\Sigma}, \Sigma \xrightarrow{f} X | E \rightarrow X,
\nabla_{E}>$ in $\mathbb{C}/\mathbb{Z}$.

\textit{i})  Since the Todd form is multiplicative for the direct sum of
  bundles with connections, if $\Sigma \rightarrow X$ is multiplied by
  $V \rightarrow pt$ then
\[
< \nabla, \Sigma \times V \rightarrow X | E_{X}, \nabla_{E} > =
\mathrm{Todd} V < \nabla_{\Sigma}, V \rightarrow X | E_{X}, \nabla>. 
\]

\textit{ii}) Fixing $E_{X},\nabla_{E}$ and varying the cycle the function
  $\phi(\nabla, V \rightarrow X) = <\nabla, V \rightarrow X | E_{X}, \nabla>$ satisfies:  there is a closed form $C$ on $X$ so that
  whenever $V \xrightarrow{f} X = \partial(W \xrightarrow{F} X)$ then
  $\phi(\nabla, V \xrightarrow{f} X) = \int_{W} F^{\ast}C \;
  \textrm{Todd} \; W \; (\textrm{mod one})$.

\underline{\textbf{Remark 3}}:  It follows that the closed form $C$ of
property \textit{ii)} has integral periods in the sense that for every
closed $W \xrightarrow{F} X$ 
\[
\int_{W} F^{\ast}C \; \textrm{Todd} \; W \hspace{1cm} \textrm{belongs to} \; \mathbb{Z},
\]
and that $C$ is unique given the values in $\mathbb{C}/\mathbb{Z}$ for SAC cycles. (We denote such forms $C$ in Section 3 by $\wedge^{even}_{integrality}$.)

\textit{iii})  If we fix $E_{X}$ but change the connection from $\nabla$ to
  $\nabla^{\prime}$, then if $CS(\nabla,\nabla^{\prime})$ denotes the
  Chern Simons difference form so that $d CS(\nabla,\nabla^{\prime}) =
  ch \nabla - ch \nabla^{\prime}$, then 
\[
<\nabla_{\Sigma} , \Sigma \xrightarrow{f} X | E_{X} , \nabla> - <
\nabla_{E} , \Sigma \xrightarrow{f} X | E_{X} , \nabla^{\prime}> =
\int_{V} f^{\ast} CS(\nabla,\nabla^{\prime}) \textrm{Todd} \; V.
\]
In particular, if $CS(\nabla, \nabla^{\prime})=$ exact, the difference is zero.

\pagebreak

\section*{\S2.  Algebraic Topology invariants in $\mathbb{Q}/\mathbb{Z}$ of complex vector bundles derived from the pairing $<\nabla_{\Sigma}, \Sigma \rightarrow X | E \rightarrow X,
  \nabla_{E}>$ in $\mathbb{C}/\mathbb{Z}$ }
Recall that the $\mathbb{Z}/2$ graded functor
  $\bar{\Omega}^{\mathbb{C}}_{\ast}(X) = \Omega^{\mathbb{C}}_{\ast}(X)
  \otimes_{\Omega^{\mathbb{C}}_{\ast}(X)}Z$ 
where the SAC bordism $\mathbb{Z}$-graded modules over SAC bordism of a point are collapsed to a $\mathbb{Z}/2$-graded functor by setting $V \cdot x = 0$ if Todd$V = 0$.  This was
introduced by Conner and Floyd
  and recall their theorem [9] implying the unexpected fact
  that $\bar{\Omega}^{\mathbb{C}}_{\ast}(\ast)$ is a
  $\mathbb{Z}/2$-graded homology theory.  Since
  $\bar{\Omega}^{\mathbb{C}}_{\ast}(X)$ is a homology theory we can
  form $\bar{\Omega}^{\mathbb{C}}_{\ast}(X,\mathbb{Z}/n)$ which can
  also be defined by SAC $\mathbb{Z}/n$-manifolds via the formula 
  $\bar{\Omega}^{\mathbb{C}}_{\ast}(X,\mathbb{Z}_{n}) =
  \Omega^{\mathbb{C}}_{\ast}(X, \mathbb{Z}/n \otimes_{\Omega_{\ast}pt}
  \mathbb{Z}$ (see below).  Then we can define the homology theory $\bar{\Omega}^{\mathbb{C}}_{\ast}(X,\mathbb{Q}/\mathbb{Z})$ as the 
\[
\underset{\stackrel{\rightarrow}{n}}{\lim} \; \bar{\Omega}^{\mathbb{C}}_{\ast}(X,\mathbb{Z}/n).  
\]

The first theorem, (whose terms are explained more fully in the Remark) is 

\underline{\textbf{Theorem AT}} (Algebraic Topology):   Complex $K$-theory is isomorphic based on $\mathbb{C}/\mathbb{Z}$ characters to Rational $\textrm{Hom} (\bar{\Omega}^{\mathbb{C}}_{\ast}(X,\mathbb{Q}/\mathbb{Z}), \mathbb{Q}/\mathbb{Z})$, whose elements are called $\mathbb{Q}/\mathbb{Z}$ characters.

\underline{\textbf{Remark}}: We will use the pairing of Section 1 
\[
<\nabla_{\Sigma}, \Sigma \rightarrow X | E \rightarrow X, \nabla_{E}> \; \in \; \mathbb{C}/\mathbb{Z}
\]
in the proof, to define $\mathbb{Z}/n$ pairings for $\Sigma$ a $\mathbb{Z}/n$ SAC bordism class of $\mathbb{Z}/n$-manifolds in $X$.
\[
<\Sigma \rightarrow X | E \rightarrow X> \; \in \mathbb{Z}/n
\]
will form an inverse limit.  This limit is uncountable, but the rationality condition will characterize the image we seek.

\underline{\textbf{Remark}}:  For a homology theory $h_{\ast}$ 
\[
C \in \textrm{Rational} \; \textrm{Hom} \; (h_{\ast}( \; , \mathbb{Q}/\mathbb{Z}),\mathbb{Q}/\mathbb{Z})
\]
means by definition the boxed commutative diagram: 

\setlength{\unitlength}{0.5cm}
\begin{picture}(24,7)\thinlines

\put(10.5,6.5){\line(1,0){10}}
\put(4,5){$\xrightarrow{\beta}$}
\put(5.5,5){$h_{\ast}(X,\mathbb{Z})$}
\put(9,5){$\rightarrow$}
\put(11,5){$h_{\ast}(X,\mathbb{Q})$}
\put(14.5,5){$\rightarrow$}
\put(16,5){$h_{\ast}(X,\mathbb{Q}/\mathbb{Z})$}
\put(21,5){$\xrightarrow{\beta}$}

\put(6,3){$\downarrow  C_{\mathbb{Z}}$}
\put(11.5,3){$\downarrow C_{\mathbb{Q}}$}
\put(17.5,3){$\downarrow C$}

\put(2.5,1){$0$}
\put(4,1){$\rightarrow$}
\put(6.5,1){$\mathbb{Z}$} 
\put(9,1){$\rightarrow$}
\put(12,1){$\mathbb{Q}$}
\put(14.5,1){$\rightarrow$}
\put(17.5,1){$\mathbb{Q}/\mathbb{Z}$} 
\put(10.5,0){\line(1,0){10}}
\put(10.5,0){\line(0,1){6.5}}
\put(20.5,0){\line(0,1){6.5}}
\put(21,1){$\rightarrow$}
\put(23,1){$0$}

\end{picture}

where the upper row is the Bockstein sequence for $h_{\ast}$.  $C$ determines $C_{\mathbb{Q}}$ uniquely when $C_{\mathbb{Q}}$ exists (see below).

Now we turn to the proof of Theorem AT using the pairing $<\nabla_{\Sigma} ; \Sigma \rightarrow X | E \rightarrow X ; \nabla_{E}> \; \in \; \mathbb{C}/\mathbb{Z}$, and the $\mathbb{Z}/n$-manifold definition of $\bar{\Omega}^{\mathbb{C}}_{\ast}(X,\mathbb{Z}/n)$.  (See Remark below for an explanation of the $\mathbb{Z}/n$-manifold definition of $\Omega^{\ast}_{\mathbb{C}}(X,\mathbb{Z}/n)$.)

\underline{\textbf{Definition 1}}:  A $\mathbb{Z}/n$-manifold is a pair $(V, \beta V)$ where boundary $V$ is the disjoint union of $n$ copies of a closed manifold $\beta V$, (read ``Bockstein of $V$'').  We say the $\mathbb{Z}/n$-manifold $(V, \beta V)$ is the boundary of $(W, \beta W)$ if $\beta V$ is the boundary of $\beta W$, and $V$ union $n$-copies of $\beta W$ glued on the $n$ boundary components $\beta V$ is a closed manifold which is the boundary of $W$.  Word picture:  a $\mathbb{Z}/n$-manifold looks like a book with $n$-pages attached along a binding $\beta V$ but whose edges are all glued together to form $V$.

\begin{figure}
\centerline{\includegraphics*[scale=0.5]{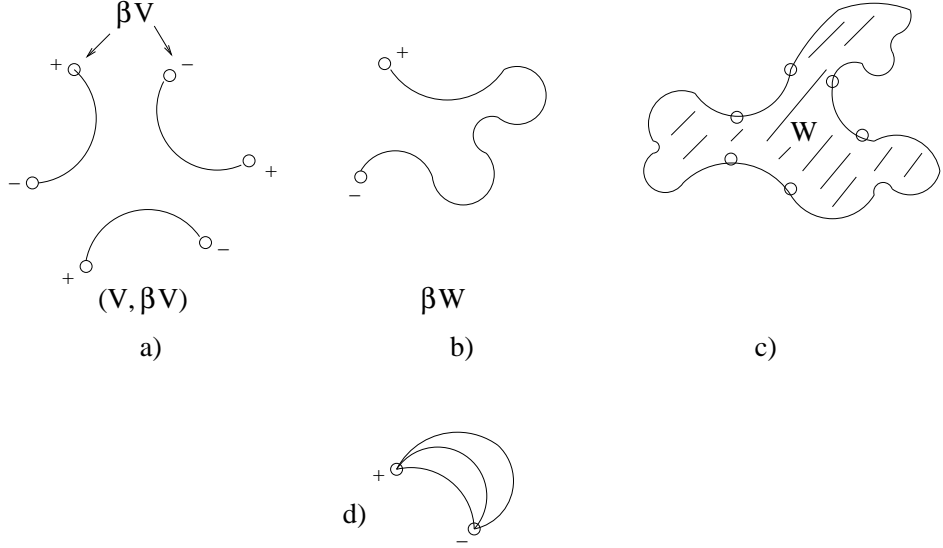}}
%\centerline{\includegraphics*[scale=0.5]{/home/ann/math4jim/paper6/p61a.eps}}
\caption{$\mathbb{Z}/n$ manifold in $X$, $n=3$}
\end{figure}

\underline{\textbf{Construction II}}:  Given an even-dimensional SAC-$\mathbb{Z}/n$ manifold in $X$, $(V, \beta V) \rightarrow X$ and a complex bundle $E \rightarrow X$, we construct an element in $\mathbb{C}/\mathbb{Z}$ of order $n$ as follows: write $\beta V = \partial Q$ and extend $E/V = E_{V}$ to $EQ$ over $Q$ using Proposition 1 of Section  1.  Enrich these objects with connections as in Section  1.  For dim $V$ even, consider the expression

$\ast) \hspace{.3cm}  \frac{1}{n} \int ch \; E_{V} \; \textrm{Todd} \; V - \int_{Q} ch \; E_{Q} \; \textrm{Todd} \; Q$

defining an element in $\mathbb{C}/\mathbb{Z}$ after reducing mod $\mathbb{Z}$.

\pagebreak

\underline{\textbf{Proposition 2}}: 

a)  The element defined in $\ast$), denoted by $<V \xrightarrow{f} X | E \rightarrow X>$ is an element of order $n$ in $\mathbb{C}/\mathbb{Z}$ and is independent of the choices $\nabla_{V}$, $\nabla_{E}$, and $Q, \nabla_{Q}$.

b)  If the $\mathbb{Z}/n$-manifold $V \rightarrow X$ in $X$ bounds as a $\mathbb{Z}/n$ manifold in $X$ then $\ast$) is zero in $\mathbb{C}/\mathbb{Z}$.

c)  $<U \rightarrow pt \times V \xrightarrow{f} X | E \rightarrow X> = \textrm{Todd} \; U < V \xrightarrow{f} X | E \rightarrow X>.$

\underline{\textbf{Proof}}:   For part a), multiply $\ast$) by $n$ and look at the boundary $R$ of Figure 1c.   This expression is the integral of $(ch \; E \; \textrm{Todd} \; R)$ where $R$ is the closed manifold ($V$ union $n$-copies of $Q$) = $R$.  This is an integer.  So $n \cdot (\ast) = 0$ in $\mathbb{C}/\mathbb{Z}$.

Changing $Q$ to another filling of $\beta V$ changes the second term in $\ast$) by an integer as in \mbox{Construction I}.  

Changing the connection can be done continuously.  Elements of order $n$ cannot move continuously.  This proves a).

To prove b), note the integer defined by [$\ast$) times $n$] is actually zero since $R = \partial W$.  Thus dividing by $n$, it is still zero as a real number.  Thus its reduction $\textrm{mod} \; \mathbb{Z}$ is zero.

c) follows from the definitions, and the multiplicative properties of the Todd form. \hspace{1 cm} $\blacksquare$

\underline{\textbf{Corollary}}:  For every $n$ we have character invariants defined using the $\mathbb{C}/\mathbb{Z}$ characters,
\[
K^{even}_{\mathbb{C}}(X) \xrightarrow{< \; | \; >} \textrm{Hom}(\bar{\Omega}_{even}^{\mathbb{C}}(X, \mathbb{Z}/n), \mathbb{Z}/n)
\]
where $\bar{\Omega}_{(\ast)}^{\mathbb{C}}$ is the Conner-Floyd homology theory ($\mathbb{Z}/2$-graded) defined by
\[
\bar{\Omega}_{(\ast)}^{\mathbb{C}}(X) = \Omega_{\ast}^{\mathbb{C}}(X) \otimes_{\Omega_{\ast}^{\mathbb{C}}pt} Z.
\]

\underline{\textbf{Proposition 3}}:  Elements on $\bar{\Omega}_{(\ast)}(X,\mathbb{Z}/n)$ have order $n$.  $\bar{\Omega}_{(\ast)}(X,\mathbb{Z}/n)$ is a multiplicative theory, (which means one can multiply cycles in $X$ and in $Y$ to get cycles in $X \times Y$).

\underline{\textbf{Proof}}:  We prove the second statement first for $\Omega_{(\ast)}^{\mathbb{C}}(X,\mathbb{Z}/n)$ (without the ``bar'').  The cartesian product of two $\mathbb{Z}/n$-manifolds (like the product of two smooth manifolds with boundary) has a codimension two locus $L$ that needs attention.  The neighborhood of the locus $L$ has the form $L \times c_{n} \times c_{n}$ where $c_{n}$ is the cone on $n$ points (denoted $(n)$).  Now $c_n \times c_n$ is the cone on the join $(n) \ast (n)$.  But $(n) \ast (n)$ defines an element in $\Omega_{1}^{\mathbb{C}}(pt,\mathbb{Z}/n)$ which is zero as seen by the exact sequence
\[ \begin{array}{ccccccccc}
(\Omega_{1}^{\mathbb{C}}(pt) & \xrightarrow{n} & \Omega_{1}^{\mathbb{C}}(pt) & \rightarrow & \Omega_{1}^{\mathbb{C}}(pt, \mathbb{Z}/n) & \rightarrow & \Omega_{0}^{\mathbb{C}} & \xrightarrow{n} & \Omega_{0}^{\mathbb{C}} ) = \\
&&&&&&&& \\
(0 & \rightarrow & 0 & \rightarrow & \Omega_{1}^{\mathbb{C}}(pt,\mathbb{Z}/n) & \rightarrow & \mathbb{Z} & \xrightarrow{n} & \mathbb{Z}). \\
\end{array} \]

\underline{\textbf{Remark}}:  For a general theory, elements in $h_{1}(pt)$ can create a difficulty at this point of the argument for Proposition 3, e.g. $n=2$ and $KO_{\ast}$. [Communication of Luke Hodgkin].

Continuing, choose a two dimensional $\mathbb{Z}/n$ SAC manifold $C_{n}$ with boundary, whose boundary is $n \ast n$.  Use it to repair the neighborhood of the locus $L$ as follows:  remove $L \times c_{n} \times c_{n}$ and replace it by $L \times C_{n}$ glued along the boundary = $L \times ((n) \ast (n))$.

Map the repaired cartesian product of cycles to the cartesian product of cycles by projecting $C_{n}$ to the cone on its boundary.  This defines (by repairing bordisms likewise) the multiplicative structure: a map, $\Omega_{\ast}^{\mathbb{C}}(X,\mathbb{Z}/n) \otimes \Omega_{\ast}^{\mathbb{C}}(Y,\mathbb{Z}/n) \rightarrow \Omega_{\ast}^{\mathbb{C}}(X \times Y,\mathbb{Z}/n)$.

Since the zero manifold [$n$ points] bounds in $\mathbb{Z}/n$ bordism of a point, that $\Omega_{\ast}^{\mathbb{C}}(X,\mathbb{Z}/n)$ is a $\mathbb{Z}/n$ module follows from the map defining the multiplicative structure.  This completes the first part of the proof of Proposition 3.  The rest of the proof identifying the $\mathbb{Z}/n$ Conner-Floyd theory with that defined by $\mathbb{Z}/n$-manifolds is in the Appendix to Section 2.

\section*{\S2. Appendix (Homology Theory)}

Continuing the proof of Proposition 3:

Let $M_{n}$ denote the $\mathbb{Z}/n$-Moore space, the circle with one two cell attached by degree $n$.  If $h_{\ast}$ is a homology theory, then $h_{k}(X,\mathbb{Z}/n)$ for $X$ connected may be defined as
\[
h_{k+1}(X \wedge M_{n}, \mathbb{Z}) \hspace{.5cm} \textrm{where} \; X \wedge M_{n} \equiv X \times M_{n} / X \vee M_{n}.
\]
Note the Bockstein exact sequence
\[ \begin{array}{ccccccc}
\xrightarrow{\beta} & h_{\ast}X & \xrightarrow{n} & h_{\ast}X & \rightarrow & h_{\ast}(X,\mathbb{Z}/n) & \xrightarrow{\beta} \\ 
\end{array} \]
follows by applying $h_{\ast}$ to the cofibration
\[
X \wedge S \xrightarrow{1 \wedge n} X \wedge S^{1} \rightarrow X \wedge M_{n}
\]
obtaining the long exact sequence of a cofibration.

That $\Omega_{\ast}^{\mathbb{C}}(X, \mathbb{Z}/n)$ defined by SAC $\mathbb{Z}/n$-manifolds agrees with this definition is proved by constructing this exact sequence directly for $\mathbb{Z}/n$-manifolds getting a map and using the 5-lemma.  This completes the discussion for Proposition 3.  

Now we can prove Proposition 4:

\underline{\textbf{Proposition 4}}:  For the homology theory,
\[
\bar{\Omega}_{\ast}^{\mathbb{C}}(X) \equiv \Omega_{\ast}^{\mathbb{C}}(X) \otimes_{\Omega_{\ast}^{\mathbb{C}}(pt)} Z,
\]
we have
\[
\bar{\Omega}_{\ast}^{\mathbb{C}}(X, \mathbb{Z}/n) = \Omega_{\ast}^{\mathbb{C}}(X,\mathbb{Z}/n) \otimes_{\Omega_{\ast}^{\mathbb{C}}(pt)} Z,
\]
where $\Omega_{\ast}^{\mathbb{C}}(X,\mathbb{Z}/n)$ is defined by SAC $\mathbb{Z}/n$-manifolds.

\underline{\textbf{Proof of Proposition 4}}:  
\[
\bar{\Omega}_{\ast}^{\mathbb{C}}(X, \mathbb{Z}/n) \equiv \bar{\Omega}_{\ast + 1}^{\mathbb{C}}(X \wedge M_{n}) \equiv \Omega_{\ast + 1}^{\mathbb{C}}(X \wedge M_{n}) \otimes_{\Omega_{\ast}^{\mathbb{C}}(pt)} Z = \Omega_{\ast + 1}^{\mathbb{C}}(X, \mathbb{Z}/n) \otimes_{\Omega_{\ast}^{\mathbb{C}}(pt)} Z 
\]
by the above.  This is one crucial place where $\bar{\Omega}$ being homology theory is used.

\underline{\textbf{Continuing Proof of Theorem AT}}:  We combine the mod $n$ character invariants of Construction II and the fact that $\bar{\Omega}_{\ast}^{\mathbb{C}}(X, \mathbb{Z}/n)$ are $\mathbb{Z}/n$ modules, i.e. each element has order $n$.

\[ \begin{array}{rcl}
K_{\mathbb{C}}^{\ast}(X,Z) & \xrightarrow{\textrm{pairing}} & \underset{\underset{n}{\leftarrow}}{\lim} \; \textrm{Hom}(\bar{\Omega}_{\ast}^{\mathbb{C}}(X, \mathbb{Z}/n),\mathbb{Z}_{n}), \textrm{ from the Corollary to Proposition 2} \\
 & = & \underset{\underset{n}{\leftarrow}}{\lim} \; \textrm{Hom}(\bar{\Omega}_{\ast}^{\mathbb{C}}(X, \mathbb{Z}/n),\mathbb{Q}/\mathbb{Z}), \textrm{ using the }\mathbb{Z}/n\textrm{-module property} \\
 & = & \textrm{Hom}(\underset{\underset{n}{\rightarrow}}{\lim} \bar{\Omega}_{\ast}^{\mathbb{C}}(X, \mathbb{Z}/n),\mathbb{Q}/\mathbb{Z}) \\
 & = & \textrm{Hom}(\Omega_{\ast}^{\mathbb{C}}(X, \mathbb{Q}/\mathbb{Z}),\mathbb{Q}/\mathbb{Z}), \\
\end{array} \]
because we have $\mathbb{Z}/n$-modules, which implies by the universal properties of the finite completion functor\footnote{$\wedge$ replaces an Abelian group by the inverse limit of its finite quotients.} $\wedge$, a commutative diagram
\[ \begin{array}{ccccc}
K_{\mathbb{C}}^{\ast}(X,Z) & \xrightarrow{\wedge} & K_{\mathbb{C}}^{\ast}(X,\hat{Z}) & \xrightarrow{<|>^{\wedge}} &  \textrm{Hom} (\bar{\Omega}_{\ast}^{\mathbb{C}} X, \mathbb{Q}/\mathbb{Z}, \mathbb{Q}/\mathbb{Z} ) \\
&&&& \\
\text{\footnotesize{Id}} \updownarrow &&&& \updownarrow \text{\footnotesize{Id}} \\
&&&& \\
K_{\mathbb{C}}^{\ast}(X,Z) & & \xrightarrow{<|>} & & \textrm{Hom} (\bar{\Omega}_{\ast}^{\mathbb{C}} X, \mathbb{Q}/\mathbb{Z}, \mathbb{Q}/\mathbb{Z} ). \\
\end{array} \]
$<|>^{\wedge}$ is a map of cohomology theories which on a point maps $\hat{\mathbb{Z}} \rightarrow \textrm{Hom}(\mathbb{Q}/\mathbb{Z}, \mathbb{Q}/\mathbb{Z})$ by an isomorphism.  Thus $<|>^{\wedge}$ is an isomorphism.

\underline{\textbf{Remark}}:  $\textrm{Hom}(\bar{\Omega}_{\mathbb{C}}^{\ast}(X, \mathbb{Z}/n), \mathbb{Z}_{n})$ is a cohomology theory since $\mathbb{Z}_{n}$ is an injective $\mathbb{Z}/n$ module.  The inverse limit of finite cohomology theories is also a cohomology theory.  Putting these together yields $\textrm{Hom}(\bar{\Omega}_{\mathbb{C}}^{\ast}(X, \mathbb{Q}/\mathbb{Z}), \mathbb{Q}/\mathbb{Z})$, is a cohomology theory.  A second way is $\mathbb{Q}/\mathbb{Z}$ is divisible so Hom(homology theory, $\mathbb{Q}/\mathbb{Z}$)\footnote{We referred to this construction as the ``Pontryagin dual'' cohomology theory [13] which we learned about from Don Anderson who never published it to my knowledge.  It is now referred to by specialists as the ``Anderson dual''.  It's a great construction.} satisfies exactness and thus all the axioms to be a cohomology theory.

We also have $K_{\mathbb{C}}^{\ast}X \xrightarrow{<|>_{\mathbb{Z}}} \textrm{Hom}(\bar{\Omega}_{\ast}^{\mathbb{C}},\mathbb{Z})$ defined by:  for any closed even-dimensional SAC $V \xrightarrow{f} X$, form the integer $\int_{V} ch \; E \; \textrm{Todd} V$ where $E$ is any complex bundle over $X$.  One knows $<|>_{\mathbb{Z}} \otimes Q$ is an isomorphism of $(K_{\mathbb{C}}^{\ast}X) \otimes Q$ with $\textrm{Hom} (\bar{\Omega}_{\ast}^{\mathbb{C}}(X,\mathbb{Q}),\mathbb{Q})$.

Furthermore, for any complex bundle $E$, the diagram

\setlength{\unitlength}{0.5cm}
\begin{picture}(20,7)\thinlines
\put(0,7){\line(1,0){13}}
\put(1,5){$\bar{\Omega}_{\ast}^{\mathbb{C}}(X,\mathbb{Q})$}
\put(6,5){$\rightarrow$}
\put(8,5){$\bar{\Omega}_{\ast}^{\mathbb{C}}(X,\mathbb{Q}/\mathbb{Z})$}

\put(1,3){$<|E>_{\mathbb{Q}} \downarrow$}
\put(8,3){$\downarrow <|E>_{\mathbb{Q}/\mathbb{Z}}$}
\put(14,3){= \textit{diagram}}

\put(2,1){$\mathbb{Q}$}
\put(6,1){$\rightarrow$}
\put(8,1){$\mathbb{Q}/\mathbb{Z}$} 
\put(0,0){\line(1,0){13}}
\put(0,0){\line(0,1){7}}
\put(13,0){\line(0,1){7}}

\end{picture}

commutes.

The map $K_{\mathbb{C}}^{\ast} \rightarrow \{\mathrm{diagram\}}$ is a map of cohomology theories (see next Remark) which for a point is

\[ \begin{array}{cccccc} 
& & \mathbb{Q} & \xrightarrow{\mathrm{mod} 1} & \mathbb{Q}/\mathbb{Z} & \\
n \in \mathbb{Z} & \longmapsto & \downarrow n & & \downarrow n & ,\\
& & \mathbb{Q} & \xrightarrow{\mathrm{mod} 1} & \mathbb{Q}/\mathbb{Z} & \\
\end{array} \]
an isomorphism.  Thus it is an isomorphism of functors.  This proves Theorem AT. \hspace{1 cm} $\blacksquare$

\underline{\textbf{Remark}}:  Diagrams are easily identified (using iso $<|>^{\wedge}$) with the kernel of 
\[
(K_{\mathbb{C}}^{\ast}(X) \otimes Q) \oplus K_{\mathbb{C}}^{\ast}(X,\hat{\mathbb{Z}}) \xrightarrow{\Delta} K_{\mathbb{C}}^{\ast}(X,\hat{\mathbb{Z}}) \otimes Q
\]
where $\Delta$ is the difference of the natural maps.  By the fibre construction and the exactness of 
\[ \begin{array}{ccccccccc}
0 & \rightarrow & \mathrm{kernel} & \rightarrow  & (\;) \oplus (\;) & \rightarrow & K_{\mathbb{C}}^{\ast}(X,\hat{\mathbb{Z}}) \otimes Q & \rightarrow & 0 \\
\end{array} \]
we get that diagrams form a cohomology theory.  See [14].

\section*{\S3. The maps from differential $K$-theories to $\hat{K}$-characters are bijections ($\mathbb{C}$-linear and unitary cases)}

The $\mathbb{C}$-linear case in detail:

\underline{\textbf{Definition}}:  A $\hat{K}$-character on $X$ is an additive over disjoint union function from (enriched closed odd-dimensional SACs mapping to $X$) to ($\mathbb{C}/\mathbb{Z}$) satisfying properties i) and ii) of Remark 2 following Construction I in Section 1, using $\mathbb{C}$-valued differential forms.

\underline{\textbf{Proposition 3}}:  The construction of a $\hat{K}$-character for a complex vector bundle $(E,\nabla)$ over $X$ with $\mathbb{C}$-linear connection in Section 1 only depends on the equivalence class of $(E, \nabla)$ in differential $K$-theory, $\hat{K}(X)$ defined using $\mathbb{C}$-valued differential forms.

\underline{\textbf{Proof}}:  By Remark 2 property \textit{iii}) of Section 1, if $CS(\nabla,\nabla^{\prime})$ is exact, the $\hat{K}$-characters of $(E, \nabla)$ and $(E, \nabla^{\prime})$ are equal.  Since $ch E$ is additive, we can pass to the Grothdieck group of $(E, \nabla)$ up to $CS$ equivalence, which is the definition of $\hat{K}(X)$. \hspace{1 cm} $\blacksquare$

\underline{\textbf{Theorem DG (Differential Geometry)}}:  The natural map produced by Construction I of Section 1, $\hat{K}X \rightarrow \hat{K}$-characters with values in $\mathbb{C}/\mathbb{Z}$, is an isomorphism, where $\hat{K}X$ is defined by $\mathbb{C}$-linear connections and $\mathbb{C}$-valued differential forms.

\underline{\textbf{Proof}}:  From [5] backstopped by [4] for the odd form lemma in the $\mathbb{C}$-linear case, we have the diagram
\begin{center}
\setlength{\unitlength}{0.5cm}
\begin{picture}(24,16)\thicklines

\put(5,1){$0$}
\put(20.5,1){$0$}

\put(6,2){\vector(1,1){1.5}}
\put(18,3.5){\vector(1,-1){1.5}}

\put(1,4.5){$K^{odd}_{\mathbb{C}}(X)$}
\put(7.5,4.5){$\wedge^{odd}/{\wedge^{odd}_{integrality}}$}
\put(12.5,4.75){\vector(1,0){2.5}}
\put(16.5,4.5){$\wedge^{even}_{integrality}$}
\put(13,5){\small{$d$}}

\put(2.5,6){\vector(1,1){1.5}}
\put(6.5,7.5){\vector(1,-1){1.5}}
\put(10.5,7){\small{$i_{2}$}}
\put(10.5,6){\vector(1,1){1.5}}
\put(14.5,7.5){\vector(1,-1){1.5}}
\put(15.5,7){\small{$\hat{ch}$}}
\put(18.5,6){\vector(1,1){1.5}}
\put(21.5,7.5){\vector(1,-1){1.5}}

\put(3,8){$H^{odd}(X,\mathbb{C})$}
\put(12.75,8){$\hat{K}^{even}_{(X)}$}
\put(20,8){$H^{even}(X,\mathbb{C})$}

\put(2,11){\vector(1,-1){1.5}}
\put(6,9.5){\vector(1,1){1.5}}
\put(10.5,11){\vector(1,-1){1.5}}
\put(11.5,10.5){\small{$i_{1}$}}
\put(14,10.5){\small{$\delta_{2}$}}
\put(14.0,9.5){\vector(1,1){1.5}}
\put(18,11){\vector(1,-1){1.5}}
\put(19,10.5){\small{$ch$}}
\put(21,9.5){\vector(1,1){1.5}}

\put(8,12){$\textrm{ker} \; \hat{ch}$}
\put(12,12){\vector(1,0){2.5}}
\put(15.5,12){$K^{even}_{\mathbb{C}}(X)$}

\put(5.5,14.5){\vector(1,-1){1.5}}
\put(18,13){\vector(1,1){1.5}}

\put(4.5,15){$0$}
\put(20,15){$0$}
\end{picture}
\end{center}

\pagebreak

See also remarks in the Introduction about the $\mathbb{C}$-linear connections and differential $K$-theory defined using $\mathbb{C}$-valued forms.
From Construction I of Section 2, we get natural maps of this diagram into the diagram (whose exactness will be demonstrated)
\begin{center}
\setlength{\unitlength}{0.5cm}
\begin{picture}(24,16)\thicklines

\put(5,1){$0$}
\put(20.5,1){$0$}

\put(6,2){\vector(1,1){1.5}}
\put(18,3.5){\vector(1,-1){1.5}}

\put(0,5){\small{Rational Hom}}
\put(0,4){\small{$(\bar{\Omega}_{odd}(X, \mathbb{Q}/\mathbb{Z}),\mathbb{Q}/\mathbb{Z})$}}
\put(8,4.5){$\wedge^{odd}/{\wedge^{odd}_{integrality}}$}
\put(12,5){\vector(1,0){2.5}}
\put(16.5,4.5){$\wedge^{even}_{integrality}$}
\put(13,5.5){\small{$d$}}

\put(2.5,6){\vector(1,1){1.5}}
\put(6.5,7.5){\vector(1,-1){1.5}}
\put(10.5,7){\small{$i_{2}$}}
\put(10.5,6){\vector(1,1){1.5}}
\put(14.5,7.5){\vector(1,-1){1.5}}
\put(15.5,7){\small{$\hat{ch}$}}
\put(18.5,6){\vector(1,1){1.5}}
\put(21.5,7.5){\vector(1,-1){1.5}}

\put(3,8){$\textrm{Hom}(\bar{\Omega}_{odd}X,\mathbb{C})$}
\put(11,8){$\hat{K}^{even}$-\small{characters}}
\put(18,8){$\textrm{Hom}(\bar{\Omega}_{even}X,\mathbb{C})$}

\put(2,11){\vector(1,-1){1.5}}
\put(5.25,10.5){\small{$\textrm{mod} 1$}}
\put(6,9.5){\vector(1,1){1.5}}
\put(10.5,11){\vector(1,-1){1.5}}
\put(11.5,10.5){\small{$i_{1}$}}
\put(14,10.5){\small{$\delta_{2}$}}
\put(14.0,9.5){\vector(1,1){1.5}}
\put(18,11){\vector(1,-1){1.5}}
\put(19,10.5){\small{$ch$}}
\put(21,9.5){\vector(1,1){1.5}}

\put(5,12){$\textrm{Hom}(\bar{\Omega}_{odd}^{\mathbb{C}}X,\mathbb{C}/\mathbb{Z})$}
\put(12.5,12.5){\small{$\delta_{2} \cdot i_{1}$}}
\put(12,12){\vector(1,0){2.5}}
\put(15.5,12.5){\small{Rational Hom}}
\put(15.5,11.5){\small{$(\bar{\Omega}_{even}(X, \mathbb{Q}/\mathbb{Z}),\mathbb{Q}/\mathbb{Z})$}}

\put(5.5,15){\vector(1,-1){1.5}}
\put(18,13.5){\vector(1,1){1.5}}

\put(4.5,15.5){$0$}
\put(20,15.5){$0$}
\end{picture}
\end{center}
Note:  $\wedge^{even}_{integrality}$ and $\wedge^{odd}_{integrality}$ are those closed forms $C$ whose cohomology classes satisfy integrality conditions, namely $\int_{M} \; T(M) \; C$ is an integer dim $M$ even or odd respectively.

Notational Diagram:
\begin{center}
\setlength{\unitlength}{0.5cm}
\begin{picture}(24,10)\thicklines

\put(2,0.5){8}
\put(2.15,0.75){\circle{1}}
\put(2.15,0.75){\circle{1.3}}
\put(9,0.5){7}
\put(9.2,0.75){\circle{1.3}}

\put(12,1){\vector(1,0){2.5}}
\put(17,0.5){2}
\put(17.15,0.75){\circle{1.3}}

\put(2.5,2){\vector(1,1){1.5}}
\put(6.5,3.5){\vector(1,-1){1.5}}
\put(10.5,2){\vector(1,1){1.5}}
\put(14.5,3.5){\vector(1,-1){1.5}}
\put(18,2){\vector(1,1){1.5}}

\put(5,4){6}
\put(5.2,4.25){\circle{1}}
\put(5.2,4.25){\circle{1.3}}
\put(13,4){1}
\put(13.15,4.25){\circle{1.3}}
\put(20.2,4){5}
\put(20.4,4.25){\circle{1}}
\put(20.4,4.25){\circle{1.3}}

\put(6,5.5){\vector(1,1){1.5}}
\put(10.5,7){\vector(1,-1){1.5}}
\put(14.0,5.5){\vector(1,1){1.5}}
\put(18,7){\vector(1,-1){1.5}}

\put(9,8){3}
\put(9.2,8.25){\circle{1.3}}
\put(12,8){\vector(1,0){2.5}}
\put(17,8){4}
\put(17.2,8.25){\circle{1}}
\put(17.2,8.25){\circle{1.3}}

\end{picture}
\end{center}

Note to Reader:  The double circles mean we enter the proof knowing we have isomorphisms at these locations.  We have to fight hard for positions 1 and 3, but only a little for positions 7 and 2.

Here $\bar{\Omega}_{\ast}^{\mathbb{C}}$ denotes the $\mathbb{Z}/2$-graded functor introduced by Connor and Floyd
\[
\bar{\Omega}_{\ast}^{\mathbb{C}}(X,Z) = \Omega_{\ast}^{\mathbb{C}}(X,Z) \otimes_{\Omega_{\ast}^{\mathbb{C}}(pt)} Z = \Omega_{\ast}^{\mathbb{C}}(X,Z) / (V \cdot x - \textrm{Todd} \; V \cdot x) 
\]
where the ring of integers $Z$ is a $\Omega_{\ast}^{\mathbb{C}}(pt)$ module via the ring homomorphism $\Omega_{\ast}^{\mathbb{C}}(pt) \rightarrow Z$ provided by the Todd genus and ``Rational Hom'' is discussed again momentarily.  Note: we also write $\bar{\Omega}_{\ast}^{\mathbb{C}}(X)$ for $\bar{\Omega}_{\ast}^{\mathbb{C}}(X,Z)$, and $\bar{\Omega}_{\ast}^{\mathbb{C}}(X, \mathbb{Q}/Z)$ for $\Omega_{\ast}^{\mathbb{C}}(X, \mathbb{Q}/Z) \otimes_{\Omega_{\ast}^{\mathbb{C}}(pt)} Z$.

The proof of Theorem DG depends on

\underline{\textbf{Theorem 2}}:  For $X$ a finite complex $K_{\mathbb{C}}^{even}(X)$ is a finitely generated Abelian group whose torsion is identified by Constructions I and II to $\textrm{Hom}(\mathrm{torsion} \; \bar{\Omega}_{odd}^{\mathbb{C}}(X,Z), \mathbb{C}/\mathbb{Z})$ and whose quotient by torsion is identified to $\textrm{Hom}(\bar{\Omega}_{even}^{\mathbb{C}}(X,Z),\mathbb{Z}).$  Also the same statements hold reversing even and odd.

\underline{\textbf{Proof}}:  In the Appendix to Section 2, one shows using $Z/n$-manifolds and Construction II that one has a bijection between $K^{\ast}(X)$ and Rational Hom$(\bar{\Omega}_{\ast}(X, \mathbb{Q}/\mathbb{Z}), \mathbb{Q}/\mathbb{Z})$ for $\ast$ even or odd.  Recall an element $C$ in Rational Hom$(\;,\;)$ is an element $C$ in Hom which is part of a commutative diagram 

\setlength{\unitlength}{0.5cm}
\begin{picture}(27,7)\thinlines
\put(10.75,6.5){\line(1,0){12.25}}
\put(4,5){$\xrightarrow{\beta}$}
\put(5.5,5){$\bar{\Omega}_{even}^{\mathbb{C}}(X)$}
\put(9,5){$\xrightarrow{\otimes\mathbb{Q}}$}
\put(11,5){$\bar{\Omega}_{even}^{\mathbb{C}}(X) \otimes \mathbb{Q}$}
\put(16.5,5){$\rightarrow$}
\put(18,5){$\bar{\Omega}_{even}^{\mathbb{C}}(X,\mathbb{Q}/\mathbb{Z})$}
\put(23.5,5){$\xrightarrow{\beta}$}
\put(25,5){$\bar{\Omega}_{odd}^{\mathbb{C}}(X,Z)$}

\put(6,3){$\downarrow  C_{\mathbb{Z}}$}
\put(12.5,3){$\downarrow C_{\mathbb{Q}}$}
\put(19.5,3){$\downarrow C$}

\put(2.5,1){$0$}
\put(4,1){$\rightarrow$}
\put(6.5,1){$\mathbb{Z}$} 
\put(9.5,1){$\rightarrow$}
\put(13,1){$\mathbb{Q}$}
\put(16.5,1){$\rightarrow$}
\put(19.5,1){$\mathbb{Q}/\mathbb{Z}$} 
\put(10.75,0){\line(1,0){12.25}}
\put(10.75,0){\line(0,1){6.5}}
\put(23,0){\line(0,1){6.5}}
\put(23.5,1){$\rightarrow$}
\put(25,1){$0$.}

\end{picture}

(Similarly interchanging even and odd.)

If $C_{\mathbb{Q}}$ exists given $C$, $C_{\mathbb{Q}}$ must be unique.  This follows since the difference of two would map to zero in $\mathbb{Q}/\mathbb{Z}$ so it would factor through $\mathbb{Z} \subset \mathbb{Q}$ which is impossible since the domain is a $\mathbb{Q}$ vector space.

Also the unique $C_{\mathbb{Q}}$ that fits with $C$ determines $C_{\mathbb{Z}}$.  Thus $C$ determines $C_{\mathbb{Q}}$ and $C_{\mathbb{Q}}$ determines $C_{\mathbb{Z}}$.

Conversely given any $C_{\mathbb{Z}}$ it determines $C_{\mathbb{Q}}$ by $C_{\mathbb{Q}} = C_{\mathbb{Z}} \otimes C_{\mathbb{Q}}$, which in turn determines $C$ partially on kernel $\beta$.  Since $\mathbb{Q}/\mathbb{Z}$ is divisible (and thus an ``injective $\mathbb{Z}$-module'') any such partial $C$ extends (non-uniquely) to a full $C$.  This proves the second part of Theorem 2.

Note if given $C$, $C_{\mathbb{Q}}$ were zero, then $C$ factors through image $\beta$ which is the torsion of $\bar{\Omega}_{odd}^{\mathbb{C}}(X,\mathbb{C}).$  This proves the first part of Theorem 2. \hspace{1 cm} $\blacksquare$

\underline{\textbf{Corollary}}:  A cohomology class $c$ in $H^{even}(X, \mathbb{Q})$ is the Chern character of a complex bundle over $X$ if and only for every closed even-dimensional SAC mapping to $X$, $V \xrightarrow{f} X$, $\int f^{\ast} \; c \; \textrm{Todd} \; V$ is an integer.  A similar statement holds for the transgressed $ch$ in $U$, odd-dimensional closed SACs in $X$, elements in $H^{odd}(X, \mathbb{Q})$ and maps $X \rightarrow U$.

\underline{\textbf{Proof}}: This is just unraveling the statement of the second part of Theorem 2.  The odd case follows using the suspension isomorphism $h^{\ast}(X) = h^{\ast + 1}(\Sigma X)$ applied to $K^{\ast}_{\mathbb{C}}$ and $\bar{\Omega}_{\ast}^{\mathbb{C}}(X, \mathbb{Q}/\mathbb{Z})$.  \hspace{1 cm} $\blacksquare$

\underline{\textbf{Proof of Theorem DG}}:  In the upper diagram, the diagonal sequences and the outer sequences are exact by [5].  By the Corollary to Theorem 2, the maps at positions \raisebox{.7pt}{\textcircled{\raisebox{-.9pt}{7}}} and \raisebox{.7pt}{\textcircled{\raisebox{-.9pt}{2}}} are isomorphisms.  It follows that \mbox{$0 \rightarrow \text{\raisebox{.7pt}{\textcircled{\raisebox{-.9pt}{3}}}} \rightarrow \text{\raisebox{.7pt}{\textcircled{\raisebox{-.9pt}{1}}}} \rightarrow \text{\raisebox{.7pt}{\textcircled{\raisebox{-.9pt}{2}}}} \rightarrow 0$} is exact for the middle diagram.

\underline{\textbf{Claim}}:  For the middle diagram, the upper sequence
$\text{\raisebox{.7pt}{\textcircled{\raisebox{-.9pt}{8}}}} \rightarrow
\text{\raisebox{.7pt}{\textcircled{\raisebox{-.9pt}{6}}}} \rightarrow
\text{\raisebox{.7pt}{\textcircled{\raisebox{-.9pt}{3}}}} \rightarrow
\text{\raisebox{.7pt}{\textcircled{\raisebox{-.9pt}{4}}}} \rightarrow
\text{\raisebox{.7pt}{\textcircled{\raisebox{-.9pt}{5}}}}$ is exact.  

Proof of Claim is below.

\underline{\textbf{Corollary of Claim}}:  The map is an isomorphism at position \text{\raisebox{.7pt}{\textcircled{\raisebox{-.9pt}{3}}}}.

\underline{\textbf{Proof of Corollary of Claim}}:  By the Appendix to Section 2, the map is an isomorphism at \text{\raisebox{.7pt}{\textcircled{\raisebox{-.9pt}{4}}}} and \text{\raisebox{.7pt}{\textcircled{\raisebox{-.9pt}{8}}}}.  It is an isomorphism at \text{\raisebox{.7pt}{\textcircled{\raisebox{-.9pt}{6}}}} and \text{\raisebox{.7pt}{\textcircled{\raisebox{-.9pt}{5}}}} by direct inspection.  Thus it is an isomorphism at \text{\raisebox{.7pt}{\textcircled{\raisebox{-.9pt}{3}}}} by the 5-lemma.

\underline{\textbf{Second Corollary of Claim}}:  The completion of the proof of Theorem DG.  Apply the 5-lemma to 
$0 \rightarrow 
\text{\raisebox{.7pt}{\textcircled{\raisebox{-.9pt}{3}}}} \rightarrow
\text{\raisebox{.7pt}{\textcircled{\raisebox{-.9pt}{1}}}} \rightarrow
\text{\raisebox{.7pt}{\textcircled{\raisebox{-.9pt}{2}}}} \rightarrow 0$.  QED for Theorem DG. \hspace{1 cm} $\blacksquare$

\underline{\textbf{Proof of Claim}}:  By the second part of Theorem 2 the odd case, the image of 
$\text{\raisebox{.7pt}{\textcircled{\raisebox{-.9pt}{8}}}} \rightarrow 
\text{\raisebox{.7pt}{\textcircled{\raisebox{-.9pt}{6}}}}$ 
in the middle diagram is $\textrm{Hom}(\bar{\Omega}_{odd}X,\mathbb{Z}) \subset \textrm{Hom}(\bar{\Omega}_{odd}X,\mathbb{C})$.  This is the kernel of 
$\text{\raisebox{.7pt}{\textcircled{\raisebox{-.9pt}{6}}}} \rightarrow 
\text{\raisebox{.7pt}{\textcircled{\raisebox{-.9pt}{3}}}}$.  So we have exactness at 
$\text{\raisebox{.7pt}{\textcircled{\raisebox{-.9pt}{6}}}}$.  The image of
$\text{\raisebox{.7pt}{\textcircled{\raisebox{-.9pt}{6}}}} \rightarrow 
\text{\raisebox{.7pt}{\textcircled{\raisebox{-.9pt}{3}}}}$ is the component of the identity of the locally compact Abelian group $\textrm{Hom}(\bar{\Omega}_{odd}^{\mathbb{C}} X,\mathbb{C}/\mathbb{Z})$.  The quotient by the image is isomorphic to $\textrm{Hom}(\textrm{torsion} \; \bar{\Omega}_{odd}^{\mathbb{C}} X,\mathbb{C}/\mathbb{Z})$.  This quotient injects into 
$\text{\raisebox{.7pt}{\textcircled{\raisebox{-.9pt}{4}}}}$ 
by the first part of Theorem 2.  So we have exactness at
$\text{\raisebox{.7pt}{\textcircled{\raisebox{-.9pt}{3}}}}$.  The image of
$\text{\raisebox{.7pt}{\textcircled{\raisebox{-.9pt}{3}}}} \rightarrow 
\text{\raisebox{.7pt}{\textcircled{\raisebox{-.9pt}{4}}}}$ is the torsion of
$\text{\raisebox{.7pt}{\textcircled{\raisebox{-.9pt}{4}}}}$ again by Theorem 2.  This torsion is the kernel of
$\text{\raisebox{.7pt}{\textcircled{\raisebox{-.9pt}{4}}}} \rightarrow 
\text{\raisebox{.7pt}{\textcircled{\raisebox{-.9pt}{5}}}}$ by the proof of Theorem 2.  So we have exactness at
$\text{\raisebox{.7pt}{\textcircled{\raisebox{-.9pt}{4}}}}$.

\underline{\textbf{Corollary 1 of Proof of Theorem DG}}: \linebreak kernel$(\hat{K}X \xrightarrow{ch} \wedge^{even}_{integrality})$ is isomorphic to $\textrm{Hom}(\bar{\Omega}_{odd}^{\mathbb{C}},\mathbb{C}/\mathbb{Z})$, a complex torus of dimension the sum of the odd Betti numbers of $X$.

\underline{\textbf{Corollary 2 of Proof of Theorem DG}}:  The diagonal and outer sequences of the $\hat{K}$-character diagram are exact.

\underline{\textbf{Note}}:  Replacing $\mathbb{C}/\mathbb{Z}$ by $\mathbb{R}/\mathbb{Z}$, $\mathbb{C}$-linear connections by unitary connections, and complex valued forms by real valued forms gives the proof of Theorem DG in that case mutatis mutandis.

\section*{\S4. The Wrong Way Map on $\mathbf{\hat{K}}$-characters for a Smooth SAC Family with Complex Linear Connection on the Vertical Stable Tangent Spaces together with a Horizontal Connection}

An enriched SAC cycle on the base determines an enriched SAC cycle in the total space by pullback.  The stable tangent bundle of the pullback cycle in the total space has a natural direct sum splitting and the direct sum connection, where in the base directions we use the pullback of the connection on the cycle in the base and in the vertical directions the induced complex connection.  Now we restrict attention to even-dimensional SAC fibers over the base of the family. 

\underline{\textbf{Definition of Wrong Way Map on $\mathbf{\hat{K}}$-characters}}:
Given a function $t$ on enriched SAC cycles on the total space, we get a function $b$ on enriched SAC cycles in the base by the obvious formula:
\[
b(\textrm{base cycle}) \equiv t(\textrm{pulled-back cycle}).
\]

\underline{\textbf{Proposition}}:  If $t$ satisfies the properties of a $\hat{K}$-character on the total space then $b$ satisfies the properties of a $\hat{K}$-character on the base.

\underline{\textbf{Proof}}:  The Todd form of the pullback cycle in the total space is the wedge product of (the Todd form of the cycle in the base) with (the Todd form of the vertical).  

The same will be true for the Todd form of the pullback of a SAC enriched bordism deformation of the cycle in the base.

If $C(t)$ denotes the variation form of a $\hat{K}$-character on $T$, define $C(b)$  by the integration along the fibres of the product of $C(t)$ with the vertical Todd form on the total space.

If $W$ fills in the base cycle $V$ and $\bar{W}$ is the pullback fill in of the pulled-back cycle $\bar{V}$, then the integral of $\textrm{Todd} \; \bar{W} \cdot C(t)$ over $W$ computed by integrating along the fibres is seen to be the integral $C(b)$ over $W$.  This follows since if $I$ denotes integration along the fibres and ${\scriptstyle\Pi}$ is the projection,
\[ \begin{array}{rcl}
I(\textrm{Todd} \; \bar{W} \wedge C(t)) & = & I({\scriptstyle\Pi}^{\ast} \; \textrm{Todd}\; W \wedge \textrm{Todd(vertical)} \wedge C(t)) \\
 & & \\
 & = & \textrm{Todd} \; W \wedge I(\textrm{Todd(vertical)} \wedge C(t)) \\
 & & \\
 & = & \textrm{Todd} \; W \wedge C(b).  \hspace{1 cm} \blacksquare  \\
\end{array} \]

\section*{\S5.  The Riemannian and Unitary Case and Eta Invariants of $(X,E,\nabla)$}

We will use the APS theorem [8] to compute $(V,F)$ where $F : V \rightarrow X$ is an enriched SAC cycle in $X$.  The invariant will be the eta invariant of the $\mathrm{spin}^{c}$ Dirac operator on $V$ with coefficients in $F^{\ast} E$ reduced mod one.  Now we assume the connection $\nabla$ is unitary.

First, a SAC bundle $E$ has a canonically associated complex line bundle whose first Chern class reduces mod 2 to the second Stiefel-Whitney class of $E$.  

\underline{\textbf{Proof}}:  The top exterior power of a complex vector space $U$ is canonically isomorphic to the top exterior power of $U \oplus C$.  So we have a line bundle (functorially) associated with any SAC bundle.  Call this line bundle $L$.  The first Chern class of $L$ is the first Chern class of $E$ which reduces mod 2 to the second Stiefel-Whitney class of $E$. \hspace{1 cm} $\blacksquare$

Second, applying this to the actual tangent bundle $T(V)$ of a SAC sycle, form $T(V) \otimes L$ and its complex Clifford algebra bundle associated to a metric on $T(V)$ and a $U(1)$ metric on $L$.  There is a complex bundle $S$ which is fibrewise the irreducible complex Cifford module for the Clifford algebra on $T(V) \otimes L$ well-defined up to module isomorphism.

\underline{\textbf{Proof}}:  One knows this representation fact is equivalent to having a specific $\mathrm{spin}^{c}$ lift of the $SO(n)$ structure on $T(V)$, where $\mathrm{spin}^{c}$ is the fibre product of the diagram
\[
SO(n) \xrightarrow{``w_{2}\text{''}} RP^{\infty} \hookleftarrow RP^{1} = U(1).
\]
We have just seen using $L$ that we have a homotopy commutative diagram
\[ \begin{array}{rlcl}
 & V & \xrightarrow{c_{1}(L)} & K(Z,2) = BU(1) \\
 & & & \\
T(V) & \downarrow & & \downarrow \textrm{reduction mod 2}  \\
 & & & \\
 & B_{SO(n)} & \xrightarrow{w_{2}} & K(Z_{2},2).  \\
\end{array} \]
So given a homotopy class of homotopies making it actually commutative we have a lift $V \rightarrow B_{\mathrm{spin}^{c}}$ which is a fibre product of this diagram.  Here $B_{G}$ means classifying space for $G$.

This homotopy class of homotopies comes from the (rigidly) commutative diagram of structures,

%\begin{center}
\setlength{\unitlength}{0.5cm}
\begin{picture}(24,10)\thicklines
\put(6,0.5){$B_{SO(n)}$}
\put(9,0.75){\vector(1,0){2.5}}
\put(12.25,0.5){$B_{SO}$}
\put(15,0.75){\vector(1,0){2.5}}
\put(15.75,1.25){$w_{2}$}
\put(19,0.5){$K(Z/2,2)$} 

\put(2,1.75){\small{tangent}}
\put(3.75,3.5){\vector(1,-1){1.5}}

\put(2.5,4.25){$V$}
\put(5,4.5){\vector(3,1){7}}
\put(6,5.25){$L$}

\put(7,7){\vector(1,-1){5}}
\put(10.25,4.25){\small{inclusion}}

\put(19.25,6){\vector(0,-1){3}}
\put(19.5,4.25){\small{mod 2}}

\put(2.75,6.25){SAC}
\put(3.75,5.5){\vector(1,1){1.5}}

\put(6,7.5){$B_{U}$}
\put(9,7.75){\vector(1,0){2.5}}
\put(9.75,8.25){$c_{1}$}
\put(12.25,7.5){$B_{U(1)}$}
\put(15,7.75){\vector(1,0){2.5}}
\put(15.75,8.25){$=$}
\put(19,7.5){$K(Z,2)$}
\end{picture}
%\end{center}

One can now form the $\mathrm{spin}^{c}$ Dirac operator on the complex spinors, the sections of $S$.  This operator combines the Clifford multiplications with covariant derivatives of the induced unitary connection on $S$ (see [8]).  

By the discussion in [8], one has the spectral eta invariants of this operator with coefficients in any unitary bundle when the dimension of $V$ is odd.  By the celebrated theorem in [8], this real number defined by eigenvalues of Dirac with coefficients, zeta functions thereof and analytic continuation to zero differs by an integer from the integral of $(\textrm{Todd} \; W \cdot ch \; E)$ over $W$ where $W$ is SAC, $\nabla$ on $E$ is unitary, Dirac has coefficients in $E$ and boundary $(W) = V$.

\underline{\textbf{Corollary}}:  (Eta form of Theorem DG)

The complex angle invariants in $\mathbb{C}/\mathbb{Z}$ for the complex bundle $E$ over $X$ with unitary connection lie in $\mathbb{R}/\mathbb{Z}$ and can be defined directly for odd-dimensional SAC cycles $F:V \rightarrow X$ in $X$ enriched by Levi-Civita connections on $TV$ and a unitary connection on the canonical complex line bundle $L$ over $V$ using the eta invariant reduced mod 1 of the $\mathrm{spin}^{c}$ Dirac operator with coefficients in $F^{\ast}E$ to define a bijection
\[ \begin{array}{rcl}
 \textrm{differential } K \textrm{-theory} & \stackrel{eta}{\longleftrightarrow}  & \textrm{differential } \hat{K}\textrm{-characters}. \\
\textrm{(real forms)} & & \textrm{(values in }\mathbb{R}/\mathbb{Z}) \\
\end{array}\]

Similarly, we can give an eta computation of the push forward.  

Using the result of the Appendix to Section 5 we will be able to work with both connections on the Riemannian fibration together with a unitary connection on the line bundle $L$ associated to the SAC structure in the vertical tangent bundle.

The direct sum connection by definition computes the push forward for bundles with unitary connection.  The eta invariants of the $\mathrm{spin}^{c}$ Dirac operator relative to the rescaled Levi-Civita connections on the total space converge to these push forward values as the base becomes infinitely large relative to the fibre.  

This proves the Analytic Theorem:

\underline{\textbf{Theorem AN}}:  The invariants of the push forward of $(E, \nabla)$ for $E$ SAC and $\nabla$ unitary are computed by the limits of eta invariants mod one of the rescaled Levi-Civita connections as they converge to their adiabatic limit.

\section*{Appendix to \S5 : Adiabatic Limits in Riemannian Fibrations}

\section*{1.  Two Connections}

If $\; {\scriptstyle\Pi} : F \rightarrow M$ is a fibration of riemannian manifolds
such that
${\scriptstyle\Pi}$ is a Riemannian submersion, two natural connections are present
in $T(F)$.  The first is the Riemannian connection, $\nabla^{r}$, and
the second, $\nabla^{\oplus}$, is a direct sum connection on the
vertical tangent bundle and its orthogonal complement.

Under a stretching of the base by multiplying its metric by a constant
$\lambda$, and carrying this through to the metric on $F$ so that
$\scriptstyle{\Pi}$ remains a Riemannian submersion, $\nabla^{r}$ changes to a
connection denoted by $\nabla^{\lambda r}$. $\nabla^{\oplus}$ however
remains fixed.

In this Appendix we show that $\underset{\lambda \rightarrow \infty}\lim
\nabla^{\lambda r} = \tilde{\nabla}^{r}$ is a well defined connection,
and also show that $\tilde{\nabla}^{r}$ and $\nabla^{\oplus}$ are
equivalent.  This, in the sense that the $\mathit{CS}$ terms relating the
characteristic forms of the two connections are all exact.  We refer to Section 0 for notation.

\section*{2. Riemannian Fibrations}

Let $F \xrightarrow{\Pi} M$ be a smooth fibration over a smooth
manifold, the fibers of which are Riemannian manifolds. $x \in T(F)$
will be called \textbf{vertical} if ${\scriptstyle\Pi}_{\ast}(x) = 0$.  The
collection of vertical vectors forms a sub-bundle $\mathcal{V}
\subseteq T(F)$.  Clearly $\mathcal{V} | F_{m} = T(F_{m})$.
$\mathcal{V}$ is a Riemannian vector bundle, and we shall use $<,>$ to
denote the inner product on its fibers.

Now suppose we are given $\mathcal{H} \subseteq T(F)$ a complementary
sub-bundle to $\mathcal{V}$.  I.e. $T(F) \cong \mathcal{V} \oplus
\mathcal{H}$.  Elements of $\mathcal{H}$ will be called
\textbf{horizontal}, as will vector fields on $F$ all of whose
elements are horizontal.  Clearly a vector field $J$ on $M$ may be
lifted to a unique horizontal field ${\scriptstyle\Pi}^{\ast}(J)$ on $F$.  Such
horizontal fields on $F$ will be called \textbf{special}.  The
following is well known.

\underline{\textbf{Lemma 2.1}}:  Let $H, I$ be special horizontal
fields on $F$, and $X$ a vertical field.  
\[ \begin{array}{ll}
a) & [H,X] \; \mathrm{is \; vertical}  \\
\\
b) & [H, I] = {\scriptstyle\Pi}^{\ast}([{\scriptstyle\Pi}_{\ast}
  (H),{\scriptstyle\Pi}_{\ast} (I)]) + \mathrm{vertical}  \\
\end{array} \]

\underline{\textbf{Proof}}:  Because $H$ is special, the 1-parameter
flow induced by $H$ takes fibers to fibers.  a) is the infinitesimal
version of this observation.  To see b), let $f \in C^{\infty}(M)$.
For $p \in F$, the fact that $H, I$ are special shows 
\[ \begin{array}{ll}
 & [H,I](p)({\scriptstyle\Pi}^{\ast}(f)) =  [{\scriptstyle\Pi}_{\ast}(H), {\scriptstyle\Pi}_{\ast}(I)]({\scriptstyle\Pi}(p))(f) \\
\Longrightarrow & \\
 & {\scriptstyle\Pi}_{\ast}([H,I]) =
                {\scriptstyle\Pi}_{\ast}({\scriptstyle\Pi}^{\ast}([{\scriptstyle\Pi}_{\ast}(H),{\scriptstyle\Pi}_{\ast}(I)]))
                \hspace{2 cm} \textrm{which implies b).}  \hspace{1 cm} \blacksquare \\
\end{array}
\]
The Riemannian connection on the tangent bundles of the fibers of $F$
may be extended to an inner product preserving connection,
$\nabla^{\mathcal{V}}$, on $\mathcal{V}$ over all of $F$ as follows:

Let $X,Y,Z$ be vertical vector fields on $F$, and $H$ a special
horizontal field.  Let $\nabla$ denote the Riemannian connection on
the fibers.  Set
\[ \begin{array}{llcl}
2.2) & <\nabla^{\mathcal{V}}_{X} Y, Z> & = & <\nabla_{X} Y, Z> \\
\\
     & <\nabla^{\mathcal{V}}_{H} Y, Z> & = & \frac{1}{2} \{ <[H,Y],Z>
- <[H,Z],Y> + H(<Y,Z>) \}. \\
\end{array} \]

Direct calculation shows that $\nabla^{\mathcal{V}}$ is a well defined
connection on $\mathcal{V}$, and that $\nabla^{\mathcal{V}}$ preserves
$<,>$.

Let us now suppose that $M$ itself is a Riemannian manifold.  Since
$\mathcal{H} \cong {\scriptstyle\Pi}^{\ast}(T(M))$, the metric and the Riemannian
connection on $T(M)$ induce an inner product and connection on
$\mathcal{H}$, denoted respectively by $<,>$ and
$\nabla^{\mathcal{H}}$.  Set
\begin{eqnarray*}
\nabla^{\oplus} = \nabla^{\mathcal{V}} \oplus \nabla^{\mathcal{H}}.
\end{eqnarray*}
By making $\mathcal{V}$ and $\mathcal{H}$ orthogonal, $<,>$ becomes a
positive definite inner product on $T(F)$, on which it induces the
Riemannian connection $\nabla^{r}$.  We wish to compare $\nabla^{r}$
and $\nabla^{\oplus}$, two metric preserving connections on $T(F)$.
Letting $\mathrm{Skew}(T(F))$ denote skew symmetric endomorphisms, we
set
\begin{eqnarray*}
B = \nabla^{r} - \nabla^{\oplus} \in \wedge^{1}(F,\mathrm{Skew}(T(F))).
\end{eqnarray*}
Let $X,Y,Z$ be vertical vector fields and $H,I,J$ special horizontal
fields.

\pagebreak

\underline{\textbf{Proposition 2.3}}:  Assume we are working in a
neighborhood where the inner products of all the above pairs are
constant.  Then
\[ \begin{array}{llcl}
1) & <B_{X}Y,Z> & = & 0 \\
\\
2) & <B_{H}Y,Z> & = & 0 \\
\\
3) & <B_{X}Y,H> & = & \frac{1}{2} \{ <[H,X],Y> + <X,[H,Y]> \} \\
\\
4) & <B_{H}Y,I> & = &-\frac{1}{2} <[H,I],Y> \\
\\
5) & <B_{X}I,Z> & = &-\frac{1}{2} \{ <[I,X],Z> + <X,[I,Z]> \} \\
\\
6) & <B_{H}I,Z> & = & \; \frac{1}{2} <[H,I],Z> \\
\\
7) & <B_{X}I,J> & = &-\frac{1}{2} <[I,J],X> \\
\\
8) & <B_{H}I,J> & = & 0 \\
\end{array} \]

\underline{\textbf{Proof}}:  We recall the Koszul formula for the
Riemannian connection, as applied to the case of triples of vector
fields, the pair-wise inner products of which are constant.
\[
2.4) \qquad <\nabla_{W_{1}}W_{2}, W_{3}> =
\frac{1}{2}\{<[W_1,W_2],W_3> + <[W_3,W_1],W_2> + <[W_3,W_2],W_1>  \}
\]
To show 1) we note that
\begin{eqnarray*}
<B_{X}Y,Z> = <\nabla^{r}_{X}Y,Z> - <\nabla^{\mathcal{V}}_{X}Y,Z>.
\end{eqnarray*}
By definition of $\nabla^{\mathcal{V}}$, it was the extension of the
Riemannian connection on the fibers of $F$ to all of $T(F)$.  Since
the Riemannian connection on a submanifold is simply its orthogonal
projection to the sub-tangent bundle, $<\nabla^{r}_{X}Y,Z> =
<\nabla^{\mathcal{V}}_{X}Y,Z>$.  

2) follows immediately by comparing 2.2) to 2.4) and using the fact
that inner products of our fields are constant.

3) follows from 2.4) by noting that $\nabla^{\oplus}$ preserves each
of $\mathcal{V}$ and $\mathcal{H}$, as does 4) via a) of Lemma 2.1.

5) and 6) follow from 3) and 4) respectively, using the skew symmetric
action of the values of $B$.

To show 7), we note that $\nabla^{\mathcal{H}}$ is the pull-back under
${\scriptstyle\Pi}$ of the Riemannian connection of $T(M)$, and thus, since
${\scriptstyle\Pi}_{\ast}(X) = 0$, $\nabla^{\mathcal{H}}_{X} = 0$.  The rest follows
from 2.4) and a) of Lemma 2.1.

To show 8), we use b) of Lemma 2.1, and use 2.4) on both $T(F)$ and
$T(M)$.  Together, that shows that $<\nabla^{\mathcal{H}}_{H}I,J> = <\nabla^{r}_{H}I,J>$. \hspace{1 cm} $\blacksquare$

\section*{3. The Adiabatic Connection}

We now stretch $M$ by considering the 1-parameter family of metrics
$<,>^{\lambda} = \lambda<,>$, where \mbox{$\lambda \in [1,\infty)$}.  Lifting
  this to $F$ we see for $X,Y$ vertical and $H,I$ horizontal
\[
3.1) \qquad <X,Y>^{\lambda} = <X,Y>, \; <X,H>^{\lambda} = 0, \;
<H,I>^{\lambda} = \lambda<H,I>.
\]
Since $\nabla^{\mathcal{V}}$ is independent of a metric on $M$, and
since the Riemannian connection on $T(M)$ is unchanged under a
constant conformal change of metric, $\nabla^{\oplus} =
\nabla^{\mathcal{V}} \oplus \nabla^{\mathcal{H}}$ is invariant as
$\lambda$ changes.  The Riemannian connection on $T(F)$ does change,
however, and and we denote this family of connections by $\{
\nabla^{\lambda r}\}$.

\underline{\textbf{Theorem 3.2}}:  Set $\tilde{\nabla}^{r} =
\underset{\lambda \rightarrow \infty}{\lim} \nabla^{\lambda r}$.
Then, $\tilde{\nabla}^{r}$ is well defined, and $\tilde{\nabla}^{r}$
is equivalent to $\nabla^{\oplus}$.

\underline{\textbf{Proof}}:  Let $B^{\lambda} = \nabla^{\lambda r} -
\nabla^{\oplus} \in \wedge^{1}(F,\mathrm{End}(T(F)))$.  From 3.1) and
Proposition 2.3) we see

\[ \begin{array}{llclcl}
1) & <B^{\lambda}_{X}Y,Z> & = & 0 & & \\
\\
2) & <B^{\lambda}_{H}Y,Z> & = & 0 & & \\
\\
3) & <B^{\lambda}_{X}Y,H> & = & \frac{1}{\lambda}
<B^{\lambda}_{X}Y,H>^{\lambda} & = & \frac{1}{\lambda} <B_{X}Y,H> \\
\\
4) & <B^{\lambda}_{H}Y,I> & = & \frac{1}{\lambda}
<B^{\lambda}_{H}Y,I>^{\lambda} & = & \frac{1}{\lambda} <B_{H}Y,I>\\
\\
5) & <B^{\lambda}_{X}I,Z> & = &<B_{X}I,Z> & &  \\
\\
6) & <B^{\lambda}_{H}I,Z> & = &<B_{H}I,Z> & & \\
\\
7) & <B^{\lambda}_{X}I,J> & = &\frac{1}{\lambda} <B^{\lambda}_{X}I,J>^{\lambda} & = & \frac{1}{\lambda} <B_{X}I,J> \\
\\
8) & <B^{\lambda}_{H}I,J> & = & 0 & & \\
\end{array} \]
Setting $\tilde{B} = \underset{\lambda \rightarrow \infty}{\lim}
B^{\lambda}$, we see from the above
\[
3.3) \qquad \tilde{B}_{s} | \mathcal{V} = 0 \qquad \textrm{and} \qquad
\tilde{B}_{s}(H) = (B_{s} | H)^{\mathcal{V}}
\]
where $s$ is any tangent vector to $F$, and $( \; )^{\mathcal{V}}$ means
projection into $\mathcal{V}$.  Thus
\[
3.4) \qquad \tilde{\nabla}^{r} = \nabla^{\oplus} + \tilde{B}
\]
implying $\tilde{\nabla}^{r}$ is well defined.

Let $[\tilde{B},\tilde{B}] \in \wedge^{2}(F,\mathrm{End}(T(F)))$ be
defined as usual, i.e.  $[\tilde{B},\tilde{B}](x,y) = [B_{x},B_{y}].$
By 3.3)
\[
3.5) \qquad [\tilde{B},\tilde{B}] = 0.
\]
Let $d$ denote exterior differentiation with respect to
$\nabla^{\oplus}$ of forms on $F$ taking values in $\mathrm{End}(T(F))$.
Since $\nabla^{\oplus}$ preserves $\mathcal{V}$ and $\mathcal{H}$,
3.3) shows, for any $s$,$u$ tangent to $F$
\[
3.6) \qquad d\tilde{B}_{s,u} | \mathcal{V} = 0 \qquad \textrm{and} \qquad
d\tilde{B}_{s,u}(\mathcal{H}) \subseteq \mathcal{V}.
\]
For $t \in [0,1]$, let $\gamma(t) = \nabla^{\oplus} + t\tilde{B}$, a
curve of connections joining $\nabla^{\oplus}$ to $\tilde{\nabla}^{r}$.

Let $R^{t}$ denote the curvature tensor of $\gamma(t)$.  By the usual
formula
\begin{eqnarray*}
R^{t} = R + t d\tilde{B} + t^2[\tilde{B},\tilde{B}]
\end{eqnarray*}
and by 3.5)
\[
3.7) \qquad R^{t} = R + t d\tilde{B}.
\]
Since $\frac{d}{dt}(\gamma(t)) = \tilde{B}$, following 1.2) in $\S0$,
$\mathit{TP}_{l}(\nabla^{\oplus},\tilde{\nabla}^{r})$ consists of integrals of
terms of the form
\begin{eqnarray*}
\mathrm{tr}(\tilde{B}_{s_{1}} R^{t}_{s_{2},s_{3}} \cdots R^{t}_{s_{2l-2},s_{2l-1}}).
\end{eqnarray*}
Since $R_{s_{i},s_{j}}$ preserves $\mathcal{V}$ and $\mathcal{H}$, by
3.3), 3.6) and
3.7) we see that the endomorphism inside the parentheses is either $0$
or takes $\mathcal{H} \rightarrow \mathcal{V}$ and $\mathcal{V}
\rightarrow 0$.  In either case its trace is $0$.  Thus,
\begin{eqnarray*}
\mathit{TP}_{l}(\nabla^{\oplus},\tilde{\nabla}^{r}) = 0
\end{eqnarray*}
and thus by Proposition 1.5 in Section 0, $\nabla^{\oplus} \sim
\tilde{\nabla}^{r}$.  \hspace{1 cm} $\blacksquare$

\section*{References}

\begin{enumerate}
\item Sullivan, Dennis. ``Geometric Topology :  Localization,Periodicity and Galois Symmetry, The 1970 MIT Notes''.  Springer 2005.

\item Cheeger, Jeff and Simons, James.  ``Differential Characters and Geometric Invariants''.  Notes of Stanford Conference 1973, Lecture Notes in Math. No. 1167.  Springer-Verlag, New York.  1985.  pp. 50-90.

\item Sullivan, Dennis. ``On the Hauptvermutung for Manifolds''. Bull. Amer. Math. Soc. 73. 1967. pp. 598-600.

\item Pingali, V. and Takhtajan, Leon.  ``On Bott-Chern Forms and their Applications''.  Mathematische Annalen 360.  2014. pp. 519-546.

\item Simons, James and Sullivan, Dennis.  ``Structured Vector Bundles Define Differential $K$-theory''.  Clay Math. Proc. 11. 2010. pp. 579-599.

\item Hopkins, M.J. and Singer, I.M. ``Quadratic functions in Geometry, Topology, and $M$-theory''.  J. Diff. Geom 70. 2005. pp. 329-452.

\item Biswas, I. and Pingali, V.  ``Inverses of Structured Bundles''.  arXiv;1502.00071v2.  To appear in Geometriae Dedicata.

\item Atiyah, M.F., Patodi, V.K. and Singer, I.M. ``Spectral Asymmetry and Riemannian Geometry III''.  Math. Proc. Camb. Phil. Soc. 79.  1976. pp. 71-99.

\item Conner, Pierre and Floyd, Ed. ``The Relation of Cobordism to $K$
theory''.  1966 Springer Lecture Notes.

\item Freed D.S. and Lott, John.  ``An index theorem in differential $K$ theory''. arxiv0907.3508v.

\item Bunke, Ulrich and Schick, Thomas.  ``Uniqueness of smooth extensions of generalized cohomology theories''. arxiv submitted Jan. 29, 2009.

\item Cheeger, Jeff.  ``Eta-invariants, the adiabatic approximation
and conical singularities''.  \\
Journal of Differential Geometry 26.  1987.  pp. 175-221.

\item Sullivan, Dennis.  ``Triangulating Homotopy Equivalences and Homeomorphisms.  Geometric Topology Seminar Notes.  1967'' in ``The Hauptvermutung Book''. Ranicki, A.A., editor.  Springer.  1996.

\item Freed, D. S. ``$Z/k$ manifolds and families of Dirac operators''.  Invent. Math 92.  1988. pp. 243-254. 

\item Simons, James and Sullivan, Dennis.  ``Axiomatic Characterization of Ordinary Differential Cohomology''. arxiv:math.0701077v1.  Journal of Topology 1 no.1.  2007. 

\item Freed, D.S.  ``On Determinant Line Bundles''.  Mathematical Aspects of String Theory. ed. S.T. Yau.  World Scientific.  1987. pp. 189-238.

\item Bismut, J .M. and D.S. Freed.  ``The Analysis of Elliptic Familes II''.  Comm. Math. Phys. 107.  1986.  pp. 103-163.

\end{enumerate}

\end{document}